\documentclass{amsart}
\usepackage{graphicx}
\usepackage{amssymb}
\usepackage{amsfonts}
\swapnumbers
\newtheorem{theorem}{Theorem}[section]

\newtheorem{corollary}[theorem]{Corollary}

\theoremstyle{definition}

\newtheorem{assumption}[theorem]{Assumption}
\newtheorem{remark}[theorem]{Remark}

\numberwithin{equation}{section}
 \theoremstyle{plain}
 
 \numberwithin{equation}{section} 
 \numberwithin{figure}{section} 
 \theoremstyle{plain}
 \theoremstyle{remark}
 \newtheorem*{acknowledgement*}{Acknowledgement}

\newcommand{\cA}{{\mathcal A}}

\newcommand{\cC}{{\mathcal C}}

\newcommand{\cF}{{\mathcal F}}
\newcommand{\cG}{{\mathcal G}}
\newcommand{\cH}{{\mathcal H}}

\newcommand{\cL}{{\mathcal L}}

\newcommand{\cN}{{\mathcal N}}

\newcommand{\cX}{{\mathcal X}}

\newcommand{\Om}{{\Omega}}
\newcommand{\om}{{\omega}}
\newcommand{\ve}{{\varepsilon}}
\newcommand{\del}{{\delta}}
\newcommand{\Del}{{\Delta}}
\newcommand{\gam}{{\gamma}}
\newcommand{\Gam}{{\Gamma}}

\newcommand{\Sig}{{\Sigma}}
\newcommand{\sig}{{\sigma}}
\newcommand{\al}{{\alpha}}
\newcommand{\be}{{\beta}}
\newcommand{\ka}{{\kappa}}
\newcommand{\la}{{\lambda}}
\newcommand{\La}{{\Lambda}}

\newcommand{\up}{{\upsilon}}
\newcommand{\vrho}{{\varrho}}

\newcommand{\bbN}{{\mathbb N}}

\newcommand{\bbS}{{\mathbb S}}

\newcommand{\bbI}{{\mathbb I}}

\newcommand{\bfY}{{\bf Y}}
\newcommand{\bfX}{{\bf X}}
\newcommand{\bfZ}{{\bf Z}}

\begin{document}
\title[]{Limit theorems for numbers of multiple returns\\
in nonconventional arrays}%
 \vskip 0.1cm
 \author{ Yuri Kifer\\
\vskip 0.1cm
 Institute  of Mathematics\\
Hebrew University\\
Jerusalem, Israel}%
\address{
Institute of Mathematics, The Hebrew University, Jerusalem 91904, Israel}
\email{ kifer@math.huji.ac.il}%

\thanks{ }
\dedicatory{ }
 \date{\today}
\begin{abstract}\noindent
For a $\psi$-mixing process $\xi_0,\xi_1,\xi_2,...$ we consider the number
$\cN_N$ of multiple returns $\{\xi_{q_{i,N}(n)}\in\Gam_N,\, i=1,...,\ell\}$ to a set $\Gam_N$
for $n$ until either a fixed number $N$ or until the moment $\tau_N$ when another multiple
return $\{\xi_{q_{i,N}(n)}\in\Del_N,\, i=1,...,\ell\}$ takes place for the first time where
$\Gam_N\cap\Del_N=\emptyset$ and $q_{i,N},\, i=1,...,\ell$ are certain functions of $n$ taking
on nonnegative integer values when $n$ runs from 0 to $N$. The dependence of $q_{i,N}(n)$'s on both
$n$ and $N$ is the main novelty of the paper. Under some restrictions on the functions $q_{i,N}$ we
obtain Poisson distributions limits of $\cN_N$ when counting is until $N$ as $N\to\infty$ and geometric distributions
limits when counting is until $\tau_N$ as $N\to\infty$. We obtain also similar results in the dynamical systems setup
 considering a $\psi$-mixing shift $T$ on a sequence space $\Om$ and studying the number of multiple
 returns $\{ T^{q_{i,N}(n)}\om\in A^a_n,\, i=1,...,\ell\}$ until the first occurrence of another multiple
 return $\{ T^{q_{i,N}(n)}\om\in A^b_m,\, i=1,...,\ell\}$ where $A^a_n,\, A_m^b$ are cylinder sets of length
  $n$ and $m$ constructed by sequences $a,b\in\Om$, respectively, and chosen so that their probabilities
  have the same order.
\end{abstract}
\maketitle
\markboth{Yu.Kifer}{Nonconventional arrays}
\renewcommand{\theequation}{\arabic{section}.\arabic{equation}}
\pagenumbering{arabic}

\section{Introduction}\label{sec1}\setcounter{equation}{0}

In \cite{Ki2} we considered nonconventional arrays of the form
\begin{equation*}
S_N=\sum_{n=1}^N\prod_{j=1}^\ell T^{q_{j,N}(n)}f_j
\end{equation*}
where $T$ is a measure preserving transformation on a probability space $(\Om,\cF,P)$, $f_j$'s are bounded
measurable functions and $q_{j,N}(n)=p_jn+q_jN,\, j=1,...,\ell$ are nonnegative functions with integer $p_j$ and
$q_j$'s. It was shown there that when $p_j$'s are distinct and $T$ is weakly mixing then $\frac 1NS_N\to
\prod_{1\leq j\leq\ell}\int f_jdP$ as $N\to\infty$ in $L^2$. Without the weak mixing assumption this convergence
fails, in general, even along sufficiently dense subsequences as the following example due to Frantzikinakis shows.
Namely, take $\ell=2,\, q_{1,N}(n)=n+N,\, q_{2,N}(n)=2n,\, f_1(x)=e^{4\pi ix},\, f_2(x)=e^{-2\pi ix}$ and $T:\bbS^1\to
\bbS^1$ having the form $Te^{2\pi iy}=e^{2\pi i(y+\al)}$ for an irrational $\al$. Then $\frac 1NS_N=e^{2\pi i(x+2N\al)}$ and since $e^{4\pi iN\al}$ visits every arc on $\bbS^1$ with the frequency proportional to its length there is no convergence
of $\frac 1NS_N$ as $N\to\infty$ even along sequences having positive upper density.

The above convergence under weak mixing was a part of the proof in \cite{Ki2} of an extension of the Szemer\' edi theorem
which says that for any subset $\La$ of nonnegative integers with a positive upper density there exists $\ve>0$
and an infinite set $\cN_\La$ of positive integers with uniformly bounded gaps such that for any $N\in\cN_\La$
the interval $[0,N]$ contains not less than $\ve N$ integers $n$ with the property that $a_n+p_jn+q_jN\in\La$
for some $a_n$ and all $j=0,1,...,\ell$.

An example in \cite{Ki2} showed that when $q_{j,N}(n)$'s depend polynomially on $n$ and $N$ then weak mixing of $T$ does
not suffice for convergence of $\frac 1NS_N$ as $N\to\infty$ even when $\ell=1$ and $q_{1,N}(n)=nN$. On the other hand,
it was shown there that if we assume strong $2\ell$-mixing of $T$ then $\frac 1NS_N$ does converge in $L^2$ to
$\prod_{j=1}^\ell f_j$ provided the polynomials $q_{1,N},...,q_{\ell,N}$ are essentially distinct (as polynomials in
$n$ and $N$) and depend nontrivially on $n$.

These results motivated the study of limit theorems for such nonconventional arrays, among them the strong law of large
 numbers (convergence almost everywhere and not just in $L^2$), the central limit theorem and the Poisson limit theorem
 type results taking also into account that limit theorems for (triangular) arrays is a well studied topic in probability
 (though it seems to appear rarely in dynamical systems). In Ch. 3 of \cite{HK} we derived these types of limit theorems
 for sums of the form
 \begin{equation*}
 S_N=\sum_{n=0}^NF(\xi_{q_{1,N}(n)},...,\xi_{q_{\ell,N}(n)})
 \end{equation*}
 where $\xi_m,\, m=0,1,...$ is a sequence of random variables with sufficiently weak dependence, $F$ is a sufficiently
 regular function and $q_{i,N}(n)=p_in+q_iN$ are linear. Under certain mixing conditions we derived almost sure convergence
  of $\frac 1NS_N$ as $N\to\infty$ assuming only that the integers $p_i,\, i=1,...,\ell$ are nonzero and distinct.
 On the other hand, the convergence in distribution of $N^{-1/2}(S_N-ES_N)$ to a normal random variable
  required that each difference $q_i-q_j$ must be divisible by the greatest common divisor of $p_i$ and $p_j$
  while without this condition the variance of $N^{-1/2}(S_N-ES_N)$ may not converge as $N\to\infty$. For instance, it
  was shown in \cite{HK} that if $\xi_m,\, m=0,1,...$ is an i.i.d. sequence of random variables, $S_N=\sum_{n=1}^NF(\xi_{2n+N},\xi_{2N-2n})$
   with a symmetric function $F$ such that $EF(\xi_0,\xi_1)=0$ and
  $EF^2(\xi_0,\xi_1)>0$, then $\lim_{N\to\infty,\, N\, odd}\frac 1NS^2_N\ne \lim_{N\to\infty,\, N\, even}\frac 1NS^2_N$,
  and so there is no limiting variance which means that the central limit theorem (in the standard form) fails.

 Concerning Poisson limit theorems we considered in \cite{HK} arrays of the form
 \begin{equation*}
 S_N=\sum_{n=1}^N\prod_{j=1}^\ell\bbI_{\Gam_N}(\xi_{q_{j,N}(n)})
 \end{equation*}
 where $\xi_m,\, m=0,1,...$ is a stationary $\psi$-mixing sequence of random variables, $q_{j,N}(n)=p_jn+q_jN$ and
 $\lim_{N\to\infty}NP\{\xi_0\in\Gam_N\}=\la$. Assuming that for any nontrivial permutation $\zeta$ of $(1,2,...,\ell)$
 the matrix $(\begin{smallmatrix} p_1&p_2&...&p_\ell\\ p_{\zeta(1)}&p_{\zeta(2)}&...
&p_{\zeta(\ell)} \end{smallmatrix})$ has rank 2 we showed that $S_N$ converges in distribution to a Poisson random
variable while without the above condition this may not hold true, in general. In the dynamical systems setup we
obtained under the above condition convergence in distribution as $m\to\infty$ to Poisson random variables of
expressions having the form
\begin{equation*}
S_{N_m}=\sum_{n=1}^{N_m}\prod_{j=1}^\ell\bbI_{A_m^a}\circ T^{q_{j,N}(n)}
\end{equation*}
provided that $\lim_{m\to\infty}N_m(P(A^a_m))^\ell=\la$, where $q_{j,N}(n)=p_jn+q_jN$, $T$ is the left shift on
a sequence space $\Om$, $A^a_m$ is a cylinder of length $m$ built by a nonperiodic sequence $a\in\Om$ and $P$ is a
$\psi$-mixing $T$-invariant probability measure on $\Om$.

An extension of the strong law of large numbers to nonconventional arrays with higher degree essentially distinct
polynomials $q_{j,N}(n),\, j=1,...,\ell$ can be obtained in the same way as in Ch. 3 of \cite{HK} just relying on the
 fact that for any $k$ and $N$ the number of $n$'s satisfying each equation $q_{i,N}(n)-q_{j,N}(n)=k$ does not exceed
 the maximum of the degrees of $q_{i,N}$ and $q_{j,N}$. Another proof of this result appears in \cite{Ha}. The central
 limit theorem is more complicated and \cite{Ha} requires that in addition to the above conditions on linear $q_i$'s each
 pair of nonlinear polynomials $q_i,\, q_j$ in $n$ and $N$ with $i\ne j$ satisfies at least one of the two following
 conditions.
 The first condition says that for any $\del>0$ there exist constants $C_\del,\, N_\del>0$ and sets $\Gam_{N,\del}\subset
 [1,N]$ with cardinality not exceeding $\del N$ so that nonlinear bivariate polynomials $q_i,\, q_j$ with $i\ne j$ satisfy $inf_{N>N_\del}N^{-1}\min_{m\in[\del N,N],\, n\in[\del N,N]\setminus\Gam_{N,\del}}|q_{i,N}(m)-q_{j,N}(n)|>0$.
 This condition
 is clearly satisfied if $q_i$ and $q_j$ have different degrees as bivariate polynomials. The second condition concerns
 polynomials with the same degree, say, deg$q_i=$deg$q_j=k$, $i\ne j$. Consider the representation
  $q_{\al,N}(n)=\sum_{l=0}^k N^lQ_{\al,l}(y)$, $\al=i,j$ where $y=n/N$ and $Q_{\al,l}$'s are nonconstant polynomials
  with nonnegative integer
 coefficients whose degree do not exceed $l$. Then the condition requires that $Q_{i,k}(y)=Q_{j,k}(c_{ij}y)$ and
 $Q_{i,k-1}(y)-Q_{j,k-1}(c_{ij}y)=r_{ij}Q'_{j,k}(c_{ij}y)$ for some reals $c_{ij}>0$ and $r_{ij}$. It is clear that
 this condition is satisfied if $q_{i,N}(n)=q_i(n)$ and $q_{j,N}(n)=q_j(n)$ are univariate polynomials depending
 only on $n$.

 The present paper is devoted to two related types of limit theorems for nonconventional arrays. The
 first one is the Poisson type limit theorems as described above but with with a more general class of functions
  $q_{j,N}(n),\, j=1,...,\ell$ which include higher degree polynomials in $n$. The second type concerns limit theorems for arrays of the form
 \[
 S_N=\sum_{n=1}^{\tau_N}\prod_{j=1}^\ell\bbI_{\Gam_N}(\xi_{q_{j,N}(n)}),
 \]
 where $\tau_N=\min\{ k\geq 1:\,\prod_{j=1}^\ell\bbI_{\Del_N}(\xi_{q_{j,N}(n)})\}$ with $\Del_N,\, N\geq 1$ being another
 sequence of Borel sets. In the dynamical systems setup we have here the sums
 \[
 S_N=\sum_{k=1}^{\tau_N}\prod_{j=1}^\ell\bbI_{A^a_{n_N}}\circ T^{q_{j,N}(k)},
 \]
 where $\tau_N=\min\{ k\geq 1:\,\prod_{j=1}^\ell\bbI_{A_{m_N}^b}\circ T^{q_{j,N}(k)}=1\}$ with $b$ being another
 nonperiodic sequence and $n_N,m_N\to\infty$ as $N\to\infty$. It turns out that if $P\{\xi_0\in\Del_N\}/P\{\xi_0\in\Gam_N\}$
 and $P(A^a_{n_N})/P(A^b_{m_N})$
 are bounded away from zero and infinity then under certain conditions $S_N$ converges in distribution to a
 random variable having the geometric distribution. Observe that even for $q_{j,N}(n)$ depending just on $n$ our results
 generalize and specify somewhat both \cite{KR1} (which improved the results from \cite{Ki1}) and \cite{KR2} while the
 additional dependence on $N$ brings here additional peculiarities.

 Our results remain valid for dynamical systems possessing appropriate symbolic representations such as Axiom A
 diffeomorphisms (see \cite{Bo}), expanding transformations and some maps having symbolic representations with an
 infinite alphabet and $\psi$-mixing invariant measure such as the Gauss map with its Gauss invariant measure and
 more general $f$-expansions (see \cite{He}). A direct application of the above results in the symbolic setup yields
 the corresponding results for arrivals to elements of Markov partitions but employing additional technique (see, for
 instance, \cite{KY}) it is not difficult to extend these results for arrivals to shrinking geometric balls.

\section{Preliminaries and main results}\label{sec2}\setcounter{equation}{0}
\subsection{$\psi$-mixing processes}
Let $(\Om,\cF,P)$ be a probability space and $\cF_{mn},\, 0\leq m\leq n\leq\infty$ be a two parameter family
of $\sig$-algebras $\cF_{mn},\, 0\leq m\leq n\leq\infty$ such that $\cF_{mn}\subset\cF_{m'n'}\subset\cF$ provided
that $m'\leq m\leq n\leq n'$. Recall, that the $\psi$-dependence coefficient between two $\sig$-algebras $\cG$ and $\cH$
can be written in the form (see \cite{Br}),
\begin{eqnarray}\label{2.1}
&\psi(\cG,\cH)=\sup_{\Gam\in\cG,\Del\in\cH}\big\{\big\vert\frac {P(\Gam\cap\Del)}
{P(\Gam)P(\Del)}-1\big\vert,\, P(\Gam)P(\Del)\ne 0\big\}\\
&=\sup\{\| E(g|\cG)-E(g)\|_{L^\infty}:\, g\,\,\mbox{is}\,\,
\cH-\mbox{measurable and}\,\, E|g|\leq 1\}\nonumber
\end{eqnarray}
and the $\psi$-dependence (mixing) in the family $\cF_{mn}$ is measured by the coefficient
\begin{equation*}
\psi(n)=\sup_{m\geq 0}\psi(\cF_{0,m},\cF_{m+n,\infty}).
\end{equation*}

Our first setup includes a $\psi$-mixing identically distributed (not necessarily stationary)
sequence of random variables $\xi_0,\xi_1,...$
defined on $(\Om,\cF,P)$ which means that $\psi(1)<\infty$ and $\psi(n)\to 0$ as $n\to\infty$ where $\psi(n)$ is defined
by (\ref{2.1}) with $\cF_{mn}=\sig\{\xi_m,...,\xi_n\}$ being the minimal $\sig$-algebra for which
$\xi_m,\xi_{m+1},...,\xi_n$ are measurable. We will be counting multiple returns by the sequence
$\xi_0,\xi_1,...$ to  measurable sets $\Gam_N$ considering the sum
\begin{equation*}
 S_N=\sum_{n=1}^N\prod_{j=1}^\ell\bbI_{\Gam_N}(\xi_{q_{j,N}(n)})
 \end{equation*}
 where $q_{1,N},...,q_{\ell,N}$ are functions in $n$ taking on nonnegative integer values when $0\leq n\leq N$
 and satisfying the conditions below.
 We will be interested to show, in particular, that $S_N$ converges in distribution as $N\to\infty$ to a Poisson
 random variable provided that $\lim_{N\to\infty}NP\{\xi_0\in\Gam_N\}$ exists. A simple example from \cite{HK},
 \[
 S_{2N}=\sum_{n=0}^{2N}\bbI_{\Gam_{2N}}(\xi_n)\bbI_{\Gam_{2N}}(\xi_{2N-n})=2\sum_{n=1}^N\bbI_{\Gam_{2N}}(\xi_n)
 \bbI_{\Gam_{2N}}(\xi_{2N-n})-\bbI_{\Gam_{2N}}(\xi_N)
 \]
 shows that $S_{2N}$ can have only even limits when $\lim_{N\to\infty}NP\{\xi_0\in\Gam_N\}$ exists, and so it cannot
  converge in distribution as $N\to\infty$ to a Poisson random
  variable. Thus, certain restrictions on the polynomials $q_{1,N},...,q_{\ell,N}$ are needed.

 \begin{assumption}\label{ass2.1} $q_{1,N}(n),...,q_{\ell,N}(n)$ are functions taking on nonnegative integer values
 on integers $n,N\geq 0$, defined arbitrarily when $n>N$ and such that for some constant $K>0$ and all $N\geq 1$ the
 following properties hold true:

 (i) For any $i\ne j$, $1\leq i,j\leq\ell$ and all integers $k,l$ the numbers of integers $n,\, 0\leq n\leq N$ satisfying
 at least one of the equations
  \begin{equation}\label{2.2}
 q_{i,N}(n)-q_{j,N}(n)=k\quad{and}\quad q_{i,N}(n)=l
 \end{equation}
 do not exceed $K$;

 (ii) For any permutation $\zeta$ of $(1,2,...,\ell)$ the number of pairs $m\ne n$
  satisfying $0\leq m,n\leq N$ and solving the system of equations
 \begin{equation}\label{2.3}
 q_{i,N}(n)-q_{\zeta(i),N}(m)=0,\,\,\, i=1,...,\ell
 \end{equation}
 does not exceed $K$;

 (iii)(stronger than (ii)) The cardinality of the set $D_N$ of pairs $m\ne n$ with $0\leq m,n\leq N$ satisfying
 \[
 \max_{1\leq i\leq\ell}q_{i,N}(m)\geq\max_{1\leq i\leq\ell}q_{i,N}(n)\geq\min_{1\leq i\leq\ell}q_{i,N}(n)\geq\min_{1\leq i\leq\ell}q_{i,N}(m)
 \]
does not exceed $K$.
 \end{assumption}
 Assumption (iii) is clearly stronger than Assumption (ii) and we will need the former in the shifts setup while
 the latter will be sufficient in the $\psi$-mixing processes setup.
 Observe, next, that $\ell=2$, $q_{1,N}(n)=n$ and $q_{2,N}(n)=N-n$ in the example above do not satisfy Assumption
 \ref{ass2.1}(ii) since taking the permutation $\zeta(1)=2,\,\zeta(2)=1$ we see that the system $n-(N-m)=0$,
 $N-n-m=0$ has $N+1$ solution pairs $n,N-n$ for $n=0,1,...,N$. Note also that if $\ell=1$ then Assumption \ref{ass2.1}(ii)
  requires only that for any $N\geq 1$ there exist at most $K$ pairs $n,m$, $n\ne m$ such that $q_{1,N}(n)=q_{1,N}(m)$.
  This will hold true if, for instance, there exists $n_0\geq 1$ such that $q_{1,N}(n)$ is strictly increasing in $n$ on $[n_0,\infty)$.
  Furthermore, if $q_{i,N}(n)=r_i(n)+g_i(N),\, i=1,...,\ell$ where $r_i$'s are nonconstant, essentially distinct polynomials in $n$
 and $g_i$'s are functions of $N$, both nonnegative for $n,N\geq 0$ and taking on integer values on integers, then $q_{i,N}$'s
 satisfy Assumption \ref{ass2.1}. Indeed, the number of solutions in Assumption \ref{ass2.1}(i) is bounded by the maximal degree
 of $r_i$'s. Next, there exists an integer $n_0\geq 1$ such that all polynomials $r_i,\, i=1,...,\ell$ are strictly increasing on $[n_0,\infty)$,
  and so if $m>n\geq n_0$ then both $\max_{1\leq i\leq\ell}q_{i,N}(m)>\max_{1\leq i\leq\ell}q_{i,N}(n)$ and $\min_{1\leq i\leq\ell}q_{i,N}(m)
  >\min_{1\leq i\leq\ell}q_{i,N}(n)$ which implies Assumption \ref{ass2.1}(iii) (and so Assumption \ref{ass2.1}(ii)).

   If $q_{1,N}(n),...,q_{\ell,N}(n)$ are polynomials in $n$ then Assumption \ref{ass2.1}(i) will be satisfied provided
 $\sup_{N\geq 1}\max_{1\leq i\leq\ell}deg_nq_{i,N}<\infty$, where $deg_n$ denotes the degree of a polynomial in $n$,
 since the number of solutions in (\ref{2.2}) is bounded in this case by $\max(deg_nq_{i,N},\, deg_nq_{j,N})$.
 A sufficient condition for Assumption \ref{ass2.1}(ii) to hold true can be obtained with the help of the B\' ezout theorem
 (see \cite{Sh}, \S 2 in Ch. III or \cite{BCR}, Section 11.5) which says that if $f$ and $g$ are two nonconstant bivariate coprime
  polynomials then the system $f(x,y)=0,\, g(x,y)=0$ has no more than deg$(f)$deg$(g)$ solution pairs $x,y$. If $f$ and
  $g$ are not coprime, i.e. there exists a nonconstant polynomial $h=h(x,y)$ such that $f=h\tilde f$ and $g=h\tilde g$
  for some polynomials $\tilde f$ and $\tilde g$, then each solution of $h=0$ solves also the system $f=0,\, g=0$, and so
  the latter system may have infinitely many solutions then. Thus, if for some $N_0\geq 1$, each $N\geq N_0$ and any
  nontrivial permutation $\zeta$ of $(1,2,...,\ell)$ there exist $i\ne j$ such that the polynomials $\tilde q_{i,N}(n,m)=
  q_{i,N}(n)-q_{\zeta(i),N}(m)$ and $\tilde q_{j,N}(n,m)= q_{j,N}(n)-q_{\zeta(j),N}(m)$ are coprime and nonconstant then
  Assumption \ref{ass2.1}(ii) holds true provided we can bound uniformly in $N\geq 1$ the number of pairs $m\ne n$ which
  solve the system $q_{i,N}(n)=q_{i,N}(m)$. For this it suffices to assume that there exists $n_0\geq 1$ such that for
  all $N\geq 1$ each polynomial $q_{i,N},\, i=1,...,\ell$ is strictly increasing on $[0,\infty)$.

For any two random variables or random vectors $Y$ and $Z$ of the same
dimension denote by $\cL(Y)$ and $\cL(Z)$ their distribution and by
\[
d_{TV}(\cL(Y),\,\cL(Z))=\sup_G|\cL(Y)(G)-\cL(Z)(G)|
\]
the total variation distance between $\cL(Y)$ and $\cL(Z)$ where the supremum
is taken over all Borel sets. Our first result is the following.
\begin{theorem}\label{thm2.2} Let $\xi_0,\xi_1,\xi_2,...$ be a $\psi$-mixing sequence of
identically distributed  random variables and assume that Assumptions \ref{ass2.1}(i)-(ii) hold true.  Let $\Gam$ be a Borel set,
 $X_n=X_{n,N}=\prod_{i=1}^\ell\bbI_\Gam(\xi_{q_{i,N}(n)})$ and $S_N=S_N^\Gam=\sum_{n=1}^{N}X_{n}$.
  Then there exists a constant $C>0$ which does not depend on $\Gam$ and $N\geq 1$ and such
  that for any positive integers $M,\, N,\, R$,
 \begin{eqnarray}\label{2.4}
 &d_{TV}(\cL(S_N^\Gam),\,\mbox{Pois}(\la_N))\leq C\bigg( NR(\Phi(\Gam))^{2\ell}+NM(\Phi(\Gam))^{\ell+1}\\
 &+MR(\Phi(\Gam))^\ell +(2-(1+\psi(R))^{\ell+1})^{-2}\psi(R)\big(M^2\Phi(\Gam)+N(\Phi(\Gam))^\ell\big)\bigg),
 \nonumber\end{eqnarray}
 provided $\psi(R)<2^{\frac 1{\ell+1}}-1$, where $\Phi(\Gam)=P\{\xi_0\in\Gam\}$, $\la_N=N(\Phi(\Gam))^\ell$
  and Pois$(\la)$ denotes the Poisson distribution with the parameter $\la$.
 \end{theorem}

 \begin{corollary}\label{cor2.3} Under the conditions of Theorem \ref{thm2.2}
 suppose that in place of one set $\Gam$ we have a sequence of Borel sets $\Gam_N$ such that
 \begin{equation}\label{2.5}
 0<C^{-1}\leq N\Phi(\Gam_N)^\ell\leq C<\infty
 \end{equation}
 for some constant $C$.  Then
 \begin{equation}\label{2.6}
 d_{TV}(\cL(S_N^{\Gam_N}),\,\mbox{Pois}(\la_N))\to 0\,\,\mbox{as}\,\, N\to\infty
 \end{equation}
 where $\la_N=N(\Phi(\Gam_N))^\ell$. In particular, if
 \begin{equation}\label{2.7}
 \lim_{N\to\infty}N(\Phi(\Gam_N))^\ell=\la
 \end{equation}
 then the distribution of $S_N$ converges in total variation as $N\to\infty$
 to the Poisson distribution with the parameter $\la$.
 \end{corollary}

Now, let $\Gam$ and $\Del$ be disjoint Borel sets and set
\[
\Sig_N^{\Gam,\Del}=\sum_{k=1}^{\tau}\prod_{i=1}^\ell\bbI_\Gam(\xi_{q_{i,N}(k)})
\]
where $\tau=\tau_\Del=\min\{ n\geq 1:\,\xi_{q_{i,N}(k)}\in\Del$ for $i=1,...,\ell\}$
 writing $\tau=\infty$ if the set in braces above is empty. Denote also by Geo$(\rho),\,\rho\in(0,1)$
 the geometric distribution with the parameter $\rho$, i.e.
\[
\mbox{Geo}(\rho)\{ k\}=\rho(1-\rho)^k\,\,\mbox{for each}\,\, k\in\bbN=\{ 0,1,...\}.
\]
Now we can state
\begin{theorem}\label{thm2.4} Let $\xi_0,\,\xi_1,\,\xi_2,...$ be a $\psi$-mixing
sequence of identically distributed random variables with a marginal distribution $\Phi$ and assume that Assumptions \ref{ass2.1}(i)-(ii)
 hold true. Then for any disjoint Borel sets $\Gam,\,\Del$ with $\Phi(\Gam),\Phi(\Del)>0$
and all positive integers $M,\, N,\, R$ with $\psi(R)<2^{\frac 1{\ell+1}}-1$ we have
\begin{eqnarray}\label{2.8}
&d_{TV}(\cL(\Sig_{N}^{\Gam,\Del}),\,\mbox{Geo}(\rho))\leq C\bigg((1-(\Phi(\Del))^\ell)^N\\
&+N(\Phi(\Gam)+\Phi(\Del))^\ell((1+\psi(M))^\ell-1)+NR(\Phi(\Gam)+\Phi(\Del))^{2\ell}\nonumber\\
&+NM(\Phi(\Gam)+\Phi(\Del))^{\ell+1}+MR(\Phi(\Gam)+\Phi(\Del))^\ell+M^2(\Phi(\Gam)+\Phi(\Del))\nonumber\\
&+(2-(1+\psi(R))^{\ell+1})^{-2}\psi(R)\big( N(\Phi(\Gam)+\Phi(\Del))^\ell+M(\Phi(\Gam)+\Phi(\Del))\big)\bigg)\nonumber\\
\nonumber\end{eqnarray}
where $\rho=\frac {(\Phi(\Del))^\ell}{(\Phi(\Gam))^\ell+(\Phi(\Del))^\ell}$
and the constant $C>0$ does not depend on $\Phi(\Gam),\, \Phi(\Del)$, $M$, $N$ and $R$.
\end{theorem}

Next, let $\Gam_N,\,\Del_N,\, N=1,2,...$ be a sequence of pairs of disjoint
Borel sets such that
\begin{equation}\label{2.9}
\Phi(\Gam_N),\, \Phi(\Del_N)\to 0\,\,\mbox{as}\,\, N\to\infty\,\,\mbox{and}\,\,
0<C^{-1}\leq\frac {\Phi(\Gam_N)}{\Phi(\Del_N)}\leq C<\infty
\end{equation}
for some constant $C$.

\begin{corollary}\label{cor2.5}
Suppose that the conditions of Theorem \ref{thm2.4} concerning the
process $\xi_0,\xi_1,\xi_2,...$ and the polynomials $q_{i,N}(n),\, i=1,...,\ell$ are
satisfied. Let $\Gam_N,\,\Del_N,\, N=1,2,...$ be Borel sets satisfying
(\ref{2.9}). Then
\begin{equation}\label{2.10}
d_{TV}(\cL(\Sig_{N}^{\Gam_N,\Del_N}),\,\mbox{Geo}(\rho_N))\to 0\,\,\mbox{as}\,\,
N\to\infty
\end{equation}
where $\rho_N=(\Phi(\Gam_N))^\ell((\Phi(\Del_N))^\ell+(\Phi(\Gam_N))^\ell)^{-1}$. In
particular, if
\begin{equation}\label{2.11}
\lim_{N\to\infty}\frac {\Phi(\Del_N)}{\Phi(\Gam_N)}=\la
\end{equation}
then the distribution of $\Sig_{N}^{\Gam_N,\Del_N}$ converges in total variation
as $N\to\infty$ to the geometric distribution with the parameter
$(1+\la^\ell)^{-1}$.
\end{corollary}

\subsection{Shifts}
Our second setup consists of a finite or countable set $\mathcal{A}$,
the sequence space $\Omega=\mathcal{A}^{\mathbb{N}}$,
the $\sigma$-algebra
 $\mathcal{F}$ on $\Omega$ generated by cylinder
sets, the left shift
 $T:\Omega\rightarrow\Omega$, and a $T$-invariant probability measure
 $P$ on $(\Omega,\mathcal{F})$.
We assume that $P$ is $\psi$-mixing
with the $\psi$-dependence coefficient given by (\ref{2.1}) and (\ref{2.2})
considered with respect to the $\sig$-algebras $\cF_{mn},\, n\geq m$ generated
by the cylinder sets $\{\om=(\om_0,\om_1,...)\in\Om:\,\om_i=a_i\,$ for
$m\leq i\leq n\}$ for some $a_m,a_{m+1},...,a_n\in\cA$. Clearly, $\cF_{mn}=
T^{-m}\cF_{0,n-m}$ for $n\geq m$. For each word $a=(a_0,a_1,...,a_{n-1})\in
\cA^n$ we will use the notation $[a]=\{\om=(\om_0,\om_1,...):\, \om_i=a_i,\,
 i=0,1,...,n-1\}$ for the corresponding cylinder set. Write
 $\Omega_{P}$ for the support
of $P$, i.e.
\[
\Om_P=\{\om\in\Om:\,P[\om_0,...,\om_n]>0\,\,\mbox{for all}\,\, n\geq 0\}.
\]
 For $n\ge 1$ set $\mathcal{C}_{n}=\{[w]\::\:w\in\mathcal{A}^{n}\}$.
Since $P$ is $\psi$-mixing it follows (see \cite{KR1}, Lemma 3.1)
that there exists $\upsilon>0$ such that
\begin{equation}\label{2.12}
P(A)\le e^{-\upsilon n}\text{ for all \ensuremath{n\ge1} and
\ensuremath{A\in\mathcal{C}_{n}}}.
\end{equation}
For any $n\geq 1$ and $V\in\cF_{0,n-1}$ set
\[
\pi(V)=\min\{ k\geq 1:\, V\cap T^{-k}V\ne\emptyset\}
\]
and $S^V_N=\sum_{k=1}^N\prod_{i=1}^\ell\bbI_V\circ T^{q_{i,N}(k)}$.
Observe that always $\pi(V)\leq n$ if $V\in\cF_{0,n-1}$.

\begin{theorem}\label{thm2.6} Suppose that Assumptions \ref{ass2.1}(i) and \ref{ass2.1}(iii) are satisfied.  Then there exists a
constant $C\geq 1$ such that for any $n,\, V\in\cF_{0,n-1}$, $N$ and $R$ satisfying  $\psi(R-n)<2^{\frac 1{\ell+1}}-1$ we have,
 \begin{eqnarray}\label{2.13}
 &d_{TV}(\cL(S_N^V),\,\mbox{Pois}(\la_N))\leq C\bigg((R+n)P(V)\\
 &+\big(nP(V)+N(P(V))^{\ell}\big)\big(RP(V)+\sum_{r=\pi(V)}^{n-1}P(T^{n-r}V)\big)\nonumber\\
 &+N(P(V))^{\ell}\psi(R-n)(2-(1+\psi(R-n))^{\ell+1})^{-2}\bigg)\nonumber
 \end{eqnarray}
 where $\la_N=N(P(V))^\ell$.
\end{theorem}

 \begin{corollary}\label{cor2.7} Suppose that Assumptions \ref{ass2.1}(i) and \ref{ass2.1}(iii) are satisfied.
 Let $V_{L}\in\cF_{0,n_L-1},\, L=1,2,...$ be a sequence of sets such that
 $n_LP(V_{L})\to 0$ and $\sum_{r=\pi(V_{L})}^{n_L-1}P(T^{n_L-r}V_{L})\to 0$ as $L\to\infty$. Let $N_L\to\infty$
 as $L\to\infty$ be a sequence of integers such that $0<C^{-1}\leq\la_{L}=N_L(P(V_{L}))^\ell\leq C<\infty$ for some constant
  $C$ and all $L\geq 1$. Then
 \begin{equation}\label{2.14}
 d_{TV}(\cL(S_{N_L}^{V_{L}}),\,\mbox{Pois}(\la_{L}))\to 0\,\,\mbox{as}\,\, L\to\infty
 \end{equation}
 and if $\lim_{L\to\infty}\la_{L}=\la$ then the distribution of $S_{N_L}^{V_{N_L}}$ converges in total
 variation as $L\to\infty$ to the Poisson distribution with the parameter $\la$. In particular, if
 $V_{L}=A^\eta_{n_L}=[\eta_0,...,\eta_{n_L-1}]=\{\om\in\Om:\,\om_0=\eta_0,...,\om_{n_L-1}=\eta_{n_L-1}\}$ with $n_L\to\infty$
 as $L\to\infty$ and $\eta\in\Om_P$ is nonperiodic then $\pi(A^\eta_{n_L})\to\infty$ as $L\to\infty$ and the above
 statements hold true for such $V_{L}$'s provided the above conditions on $\la_L$ are satisfied.
 \end{corollary}

 Next, for any $V\in\cF_{0,n-1},\, V\ne\emptyset$ and $W\in\cF_{0,m-1},\, W\ne\emptyset$ define
\[
\pi(V,W)=\min\{k\geq 1:\:V\cap T^{-k}W\ne\emptyset
\mbox{ or }W\cap T^{-k}V\ne\emptyset\}\:.
\]
It is clear that $\pi(V,W)\leq m\wedge n$, and so
\[
\kappa_{V,W}=\min\{\pi(V,W),\pi(V),\pi(W)\}\leq m\wedge n
\]
where, as usual, for $n,m\ge1$ we denote $m\wedge n=\max\{m,n\}$ and $m\vee n=\min\{m,n\}$.
Set
\[
\tau_W(\omega)=\min\{k\ge 1\::\:T^{q_{i,N}(k)}\omega\in W\,\,\mbox{for}\,\, i=1,...,\ell\}
\]
with $\tau_W(\om)=\infty$ if the event in braces does not occur and define
\[
\Sig_N^{V,W}=\sum_{k=1}^{\tau_W}\prod_{i=1}^\ell\bbI_{V}\circ T^{q_{i,N}(k)}.
\]

\begin{theorem}\label{thm2.8} Assume that Assumptions \ref{ass2.1}(i) and \ref{ass2.1}(iii)
are satisfied. Then there exists a constant $C>0$ such that for any disjoint sets
 $V\in\cF_{0,n-1}$ and $W\in\cF_{0,m-1}$ with $P(V),P(W)>0$ and all integers $n,m,N,R\ge 1$ satisfying
 $\psi(R-n\vee m)<2^{\frac 1{\ell+1}}-1$ we have
\begin{eqnarray}\label{2.15}
&\quad\quad d_{TV}(\mathcal{L}\big(\Sig_N^{V,W}),Geo(\rho)\big)\le C\bigg((1-(P(W))^\ell)^N+(n\vee m)(P(V)+P(W))\\
&+N(P(V)+P(W))^\ell\big((1+\psi(n\vee m))^\ell-1+\psi(R-n\vee m)\nonumber\\
&+R(P(V)+P(W))+\sum_{r=\ka_{V,W}}^{n\vee m-1}(P(T^{n\vee m-1}V)+(P(T^{n\vee m-1}W)\big)\bigg)
\nonumber\end{eqnarray}
where $\rho=\frac{(P(W))^\ell}{(P(V))^\ell+(P(W))^\ell}$.
\end{theorem}

\begin{corollary}\label{cor2.9}
Suppose that Assumptions \ref{ass2.1}(i) and \ref{ass2.1}(iii) hold true.
Let $V_L\in\cF_{0,n_L-1}$ and $W_L\in\cF_{0,m_L-1}$, $L=1,2,...$
be two sequences of sets such that
\begin{equation}\label{2.16}
(n_L\vee m_L)(P(V_L)+P(W_L))\to 0\quad\mbox{as}\quad L\to\infty,
\end{equation}
\begin{equation}\label{2.17}
\al_L=\sum_{r=\ka_{V_L,W_L}}^{n_L\vee m_L-1}(P(T^{n_L\vee m_L-r}V_L)+P(T^{n_L\vee m_L-r}W_L))\to 0\,\,\mbox{as}\,\, L\to\infty
\end{equation}
and for some constant $C$ and all $L\geq 1$,
\begin{equation}\label{2.18}
 0<C^{-1}\leq\frac {P(V_L)}{P(W_L)}\leq C<\infty.
\end{equation}
Let $N_L,\, L=1,2,...$ be a sequence satisfying
\begin{equation}\label{2.19}
N_L(P(W_L))^\ell\to\infty\,\,\mbox{and}
\end{equation}
\begin{equation}\label{2.20}
 N_L(P(V_L)+P(W_L))^\ell(\psi(n_L\vee m_L)+P(V_L)+P(W_L)+\al_L)\to 0\,\,\mbox{as}\,\, L\to\infty.
\end{equation}
Then
\begin{equation}\label{2.21}
d_{TV}(\cL(\Sig_{N_L}^{V_L,W_L}),\, Geo(\rho_L))\to 0\,\,\mbox{as}\,\, L\to\infty
\end{equation}
where $\rho_L=(P(W_L))^\ell((P(W_L))^\ell+(P(V_L))^\ell)^{-1}$. In particular, if $\lim_{L\to\infty}\rho_L=\rho$,
then $\Sig_{N_L}^{V_L,W_L}$ converges in total variation as $L\to\infty$ to the geometric distribution with the parameter
$\rho$. Furthermore, let $V_L=A^\xi_{n_L}=[\xi_0,...,\xi_{n_L-1}]\in\cC_{n_L}$ and $W_L=A^\eta_{m_L}=[\eta_0,...,\eta_{m_L-1}]\in\cC_{m_L}$
with $n_L,m_L\to\infty$ as $L\to\infty$ and suppose that $\xi,\eta$ are not periodic and not shifts of each other. Then
\begin{equation}\label{2.22}
\ka_{A^\xi_{n_L},A^\eta_{m_L}}\to\infty\,\,\mbox{as}\,\, L\to\infty
\end{equation}
and if also
\begin{equation}\label{2.23}
n_L\wedge m_L+\ka_{A^\xi_{n_L},A^\eta_{m_L}}-n_L\vee m_L\to\infty\,\,\mbox{as}\,\, L\to\infty
\end{equation}
then (\ref{2.17}) holds true. In fact, (\ref{2.23}) is satisfied for $P\times P$-almost all $(\xi,\eta)\in\Om\times\Om$
provided
\begin{equation}\label{2.24}
2n_L\wedge m_L-n_L\vee m_L-3\up\ln(n_L\wedge m_L)\to\infty\,\,\mbox{as}\,\, L\to\infty
\end{equation}
where $\up$ is from (\ref{2.12}).
\end{corollary}

\begin{remark}\label{rem2.10}
Recall that a sequence of random variables $\xi_0,\xi_1,...$ is called $\phi$-mixing if
\begin{eqnarray*}
&\phi(n)=\sup_{m\geq 0}\big\{\big\vert\frac {P(\Gam\cap\Del)}{P(\Gam)}-P(\Del)\big\vert :\\
&P(\Gam)\ne 0,\,\Gam\in\cF_{0,m},\,\Del\in\cF_{m+n,\infty}\big\}\to\infty\,\,\mbox{as}\,\, n\to\infty
\end{eqnarray*}
where $\cF_{kl}=\sig\{\xi_k,...,\xi_l\}$. It turns out that even when $\ell=1$ (conventional setup)
and $q_{1,N}(n)=n$, in general, $\phi$-mixing does not suffice for Theorems \ref{thm2.2} and \ref{thm2.4}
and Corollaries \ref{cor2.3} and \ref{cor2.5} to hold true. Indeed, consider an i.i.d. sequence $\eta_0,\eta_1,...$
and set $\xi_{2n}=\xi_{2n+1}=\eta_n,\, n=0,1,...$. Then $\xi_0,\xi_1,...$ is a $\phi$-mixing identically distributed
sequence but, as it is easy to see, the corresponding sums in Corollaries \ref{cor2.3} and \ref{cor2.5} will converge in
distribution to random variables taking on only even integer values, and so they cannot be Poisson or geometric distributed.
Unlike the case of $\phi$-mixing identically distributed sequences of random variables discussed above,
in the shifts setup Theorems \ref{thm2.6}, \ref{thm2.8} and Corollaries \ref{cor2.7}, \ref{cor2.9} can be derived
assuming only $\phi$-mixing when $\ell=1$ by using the technique from \cite{KY}.
\end{remark}

\section{Counting arguments}\label{sec3}\setcounter{equation}{0}

Let $\cN_N$ be the set of $n\in\{ 1,...,N\}$ such that all $q_{i,N}(n),\, i=1,...,\ell$ are distinct
and set $\hat\cN_N=\{ 1,...,N\}\setminus\cN_N$, i.e.
$\hat\cN_N=\{ n\in\{ 1,...,N\}:\, q_{i,N}(n)=q_{j,N}(n)$ for some $i,j=1,...,\ell,\, i\ne j\}$.
By Assumption \ref{ass2.1}(i),
\begin{equation}\label{3.1}
\#\hat\cN_N\leq \frac 12K\ell(\ell-1)
\end{equation}
where $\#\Gam$ denotes the cardinality of a set $\Gam$.

Introduce also
\[
U_{N,M}=\{ n\in\{ 1,...,N\}:\, |q_{i,N}(n)-q_{j,N}(n)|\geq M\,\,\mbox{for all}\,\, i,j=1,...,\ell,\, i\ne j\}.
\]
By Assumption \ref{ass2.1}(i) for each pair $i\ne j$ and any $k$ there exist no more than $K$ nonnegative integers
$n$ such that $q_{i,N}(n)-q_{j,N}(n)=k$, and so
\begin{equation}\label{3.2}
\#(\{ 1,...,N\}\setminus U_{N,M})\leq KM\ell(\ell-1).
\end{equation}

We will need also the following semi-metrics between positive integers $k,l> 0$,
\[
\del_N(k,l)=\min_{1\leq i,j\leq\ell}|q_{i,N}(k)-q_{j,N}(l)|.
\]
It follows from Assumption \ref{ass2.1}(i) that for any integers $n\in\{ 1,...,N\}$ and $k\geq 0$,
\begin{equation}\label{3.3}
\#\{ m:\, \del_N(n,m)=k\}\leq K\ell(\ell-1).
\end{equation}
For any integers $M,R\geq 1$ and $1\leq n\leq N$ introduce the sets
\[
B^{M,R}_{n,N}=\{ l:\, 1\leq l\leq M,\,\del_N(l,n)<R\}\quad\mbox{and}\quad B^R_{n,N}=B^{N,R}_{n,N}.
\]
By (\ref{3.3}), for any $n$,
\begin{equation}\label{3.4}
\# B^{M,R}_{n,N}\leq\min(M,\, KR\ell(\ell-1)).
\end{equation}

Next, set
\[
\sig_N(n,m)=\max_{1\leq j\leq\ell}\min_{1\leq i\leq\ell}|q_{i,N}(n)-q_{j,N}(m)|.
\]
Then
\begin{equation}\label{3.5}
\#\{ (m,n), m\ne n:\,\sig_N(n,m)=0\,\,\mbox{and either}\,\, n\in U_{N,1}\,\,\mbox{or}\,\, m\in U_{N,1}\}\leq K.
\end{equation}
Indeed, if either $n\in U_{N,1}$ or $m\in U_{N,1}$ and
$\sig_N(n,m)=0$ then, in fact, both $n\in U_{N,1}$ and $m\in U_{N,1}$. In order to see this, suppose, for
instance, that $n\in U_{N,1}$ and $\sig_N(n,m)=0$. Then there exist permutations $\eta$ and $\zeta$
of $\{ 1,...,\ell\}$ such that $q_{\eta(1),N}(n)<q_{\eta(2),N}(n)<...<q_{\eta(\ell),N}(n)$ and
$q_{\eta(i),N}(n)=q_{\zeta(i),N}(m)$ for all $i=1,...,\ell$, and so $q_{\zeta(1),N}(m)<q_{\zeta(2),N}(m)
<...<q_{\zeta(\ell),N}(m)$. The proof is the same assuming that $m\in U_{N,1}$. Hence, $q_{i,N}(n)=
q_{\eta^{-1}\zeta(i),N}(m)$, $i=1,...,\ell$. By Assumption \ref{ass2.1}(ii) there exists no more than $K$ pairs
 $m\ne n$ solving the latter system of equations, and so (\ref{3.5}) follows.

\section{Poisson distribution limits for $\psi$-mixing processes}\label{sec4}\setcounter{equation}{0}

Set $p_{n,N}=P\{ X_{n,N}=1\}=EX_{n,N}$ and $p_{n,l,N}=P\{ X_{n,N}=1$ and $X_{l,N}=1\}=E(X_{n,N}X_{l,N})$
where $X_{n,N},\, n=1,...,N$ are the same as in Theorem \ref{thm2.2}. Then by Theorem 1 in \cite{AGG}
(warning the reader that the estimates there have an extra factor 2 due to a difference in the definition
of $d_{TV}$),
\begin{equation}\label{4.1}
d_{TV}(\cL(S_N),\,\mbox{Pois}(\la_N))\leq b_1+b_2+b_3
\end{equation}
where
\begin{equation}\label{4.2}
b_1=\sum_{n=1}^N\sum_{l\in B^R_{n,N}}p_{n,N}p_{l,N},\,\,\, b_2=\sum_{n=1}^N\sum_{n\ne l\in B^R_{n,N}}p_{n,l,N}
\end{equation}
and
\begin{equation}\label{4.3}
b_3=\sum_{n=1}^Ns_{n,N}\,\,\mbox{with}\,\, s_{n,N}=E|E(X_{n,N}-p_{n,N}|\sig\{ X_{l,N}:\, l\in\{ 1,...,N\}\setminus
B^R_{n,N}\})|.
\end{equation}

By Lemma 3.2 in \cite{KR1}, for each $n\in U_{N,M}$,
\begin{equation}\label{4.4}
p_{n,N}=P\{\xi_{q_{i,N}(n)}\in\Gam\,\,\mbox{for}\,\,i=1,...,\ell\}\leq (1+\psi(M))^\ell(\Phi(\Gam))^\ell
\end{equation}
and for any $n$, clearly,
\begin{equation}\label{4.5}
p_{n,N}\leq P\{\xi_{q_{1,N}(n)}\in\Gam\}=\Phi(\Gam).
\end{equation}
Hence, by (\ref{3.2}), (\ref{3.4}), (\ref{4.2}), (\ref{4.4}) and (\ref{4.5}),
\begin{eqnarray}\label{4.6}
&b_1\leq NKR\ell^2(1+\Phi(M))^{2\ell}(\Phi(\Gam))^{2\ell}+K^2M^2\ell^4(\Phi(\Gam))^2\\
&+KM\ell^2(N+KR\ell^2)(1+\psi(M))^{2\ell}(\Phi(\Gam))^{\ell+1}.\nonumber
\end{eqnarray}

Next, if $\del_N(n,l)=k\geq 1$ and $n,l\in U_{N,M}$ then by Lemma 3.2 in \cite{KR1},
\begin{eqnarray}\label{4.7}
&p_{n,l,N}=P\{\xi_{q_{i,N}(n)}\in\Gam\,\,\mbox{and}\,\,\xi_{q_{i,N}(l)}\in\Gam\,\,\mbox{for}\,\, i=1,...\ell\}\\
&\leq (1+\psi(M\wedge k))^{2\ell}(\Phi(\gam))^{2\ell}.\nonumber
\end{eqnarray}
If $\del_N(n,l)=k\geq 1$ and either $n\in U_{N,M}$ or $l\in U_{N,M}$ relying on Lemma 3.2 in \cite{KR1} we
see that
\begin{equation}\label{4.8}
p_{n,l,N}\leq (1+\psi(1))^{\ell+1}(\Phi(\Gam))^{\ell+1}
\end{equation}
since if, for instance, $n\in U_{N,M}$ then we have $|q_{i,N}(n)-q_{j,N}(n)|>M$ for all $i\ne j$ and , in addition,
$|q_{1,N}(l)-q_{j,N}(n)|=k\geq 1$ which yields (\ref{4.8}). If $\del_N(n,l)=k\geq 1$ and $n,l\not\in U_{N,M}$ then
$|q_{i,N}(n)-q_{j,N}(l)|=k\geq 1$ for some $i$ and $j$, and so in this case
\begin{equation}\label{4.9}
p_{n,l,N}\leq (1+\psi(k))^2(\Phi(\Gam))^2.
\end{equation}

Next, suppose that $\del_N(n,l)=0$ and either $n\in U_{N,M}$ or $l\in U_{N,M}$.
then
\begin{equation}\label{4.10}
p_{n,l,N}\leq\min(p_{n,N},\, p_{l,N})\leq (1+\psi(M))^\ell(\Phi(\Gam))^\ell.
\end{equation}
By (\ref{3.5}) there exist no more than $K$ pairs $n\ne l,\, 1\leq n,l\leq N$ such that
$\sig_N(n,l)=0$ and either $n\in U_{N,M}$ or $l\in U_{N,M}$ in which case we will rely on the estimate
 (\ref{4.10}). If, on the other hand, $\sig_N(n,l)\geq 1$ and, say, $n\in U_{N,M}$ then
 $|q_{i,N}(n)-q_{j,N}(n)|\geq M\geq 1$ for all $i\ne j$ and $|q_{i,N}(n)-q_{m,N}(l)|\geq 1$
 for all $i=1,...,\ell$ and some $1\leq m\leq\ell$. Similarly, if $\sig_N(l,n)\geq 1$ and
  $l\in U_{N,M}$ then $|q_{i,N}(l)-q_{j,N}(l)|\geq M\geq 1$ for all $i\ne j$ and
  $|q_{i,N}(l)-q_{m,N}(n)|\geq 1$ for all $i=1,...,\ell$ and some $1\leq m\leq\ell$. In both cases
  we obtain the estimate (\ref{4.8}) in view of Lemma 3.2 from \cite{KR1}.

Finally, if $\del_N(n,l)=0$ and $n,l\not\in U_{N,M}$ then by (\ref{4.5}),
\begin{equation}\label{4.11}
p_{n,l,N}\leq p_{n,N}\leq\Phi(\Gam).
\end{equation}
It follows from (\ref{3.2})--(\ref{3.4}) and (\ref{4.7})--(\ref{4.11}) that
\begin{eqnarray}\label{4.12}
&b_2\leq NKR\ell^2(1+\psi(1))^{2\ell}(\Phi(\Gam))^{2\ell}+K^2M^2\ell^4\Phi(\Gam)\\
&+KM\ell^2(N+KR\ell^2)(1+\psi(1))^{\ell+1}(\Phi(\Gam))^{\ell+1}\nonumber\\
&+K^2M^2\ell^4(1+\psi(1))^2(\Phi(\Gam))^2+2K^2MR\ell^4(1+\psi(M))^\ell(\Phi(\Gam))^\ell.
\nonumber\end{eqnarray}

Next, we claim that for any $n=1,...,N$,
  \begin{eqnarray}\label{4.13}
  &s_{n,N}\leq 2^{2(\ell+2)}(2-(1+\psi(R))^{\ell+1})^{-2}\psi(R)E|X_{n,\al}
  -p_{n,\al}|\\
  &\leq 2^{2\ell+5}(2-(1+\psi(R))^{\ell+1})^{-2}\psi(R)p_{n,N}.
  \nonumber\end{eqnarray}
  Indeed, let $\cG$ be the $\sig$-algebra generated by all $\xi_{q_{i,N}(l)},\, i=1,...,\ell$
  with $l\in\{ 1,...,N\}\setminus B^R_{n,N}$ and $\cH$ be the $\sig$-algebra
  generated by $\xi_{q_{i,N}(n)},\, i=1,...,\ell$. Since $\del_N(n,l)\geq R$
  for such $l$ and $n$ we derive from Lemma 3.3 in \cite{KR1} that
  \begin{equation}\label{4.14}
  \psi(\cG,\cH)\leq 2^{2(\ell+2)}\psi(R)(2-(1+\psi(R))^{\ell+1})^{-2}
  \end{equation}
  provided $\psi(R)<2^{\frac 1{\ell+1}}-1$ which we assume. Since
  $\sig\{ X_l:\, l\in\{ 1,...,N\}\setminus B^R_{n,N}\}\subset\cG$
  and $\sig\{ X_n\}\subset\cH$ we obtain (\ref{4.13}) from (\ref{2.1})
  and (\ref{4.14}). Now by (\ref{3.2}), (\ref{4.3})--(\ref{4.5}) and (\ref{4.13}),
  \begin{equation}\label{4.15}
  b_3\leq 2^{2\ell+5}(2-(1+\psi(R))^{\ell+1})^{-2}\psi(R)\big( KM\ell(\ell-1)\Phi(\Gam)
  +N(1+\psi(M))^\ell(\Phi(\Gam))^\ell\big).
  \end{equation}
  Finally, (\ref{4.1}), (\ref{4.6}) (\ref{4.12}) and (\ref{4.15}) yield (\ref{2.4}). \qed

  When (\ref{2.5}) holds true for $\Gam=\Gam_N$ we can choose $M=M_N\to\infty$ and $R=R_N\to\infty$
  as $N\to\infty$ so that $M_N^2\Phi(\Gam_N)\to 0$ and $M_NR_N(\Phi(\Gam_N))^\ell\to 0$ as $N\to\infty$ which
  in view of (\ref{2.4}) will yield (\ref{2.6}) proving Corollary \ref{cor2.3}.  \qed

\section{Geometric distribution limits for $\psi$-mixing processes}\label{sec5}\setcounter{equation}{0}

Set for convenience $\Gam_0=\Gam$, $\Gam_1=\Del$, $X_{n,\al}=\prod_{i=1}^\ell
\bbI_{\Gam_\al}(\xi_{q_{i,N}(n)})$, $\al=0,1$ and $S_L=\sum_{n=1}^LX_{n,0}$.
Let $X'_{n,\al},\, n=1,2,...,\, \al=0,1$ be a sequence of independent random variables such that
$X'_{n,\al}$ has the same distribution as $X_{n,\al}$. Set $\tau_N=\min(\tau,N)$, $S'_L=\sum_{n=1}^LX'_{n,0}$,
$\tau'=\min\{ n\geq 1:\, X'_{n,1}=1\}$ and $\tau'_N=\min(\tau',N)$. Next, let $Y_{n,0}$ and $Y_{n,1}$,
$n=1,2,...$ be two independent of each other sequences of i.i.d. random variables such that
\begin{equation}\label{5.1}
P\{ Y_{n,\al}=1\}=\Phi(\Gam_\al)^\ell=1-P\{ Y_{n,\al}=0\},\,\al=0,1.
\end{equation}
We can and will assume that all above random variables are defined on the
same (sufficiently large) probability space. Set also
\[
S^*_L=\sum_{n=1}^{L}Y_{n,0},\,\tau^*=\min\{ n\geq 0:\, Y_{n,1}=1\}\,\,
\mbox{and}\,\,\tau^*_N=\min(\tau^*,N).
\]

Now observe that $S^*_{\tau^*}$ has by Lemma 3.1 from \cite{KR2} the geometric
distribution with the parameter
\begin{equation}\label{5.2}
\varrho=\frac {\Phi(\Gam_1)^\ell}{\Phi(\Gam_1)^\ell+\Phi(\Gam_0)^\ell
(1-\Phi(\Gam_1)^\ell)}>\rho.
\end{equation}
Next, we can write
\begin{equation}\label{5.3}
d_{TV}(\cL(S_\tau),\,\mbox{Geo}(\rho))\leq A_1+A_2+A_3+A_4+A_5
\end{equation}
where $A_1=d_{TV}(\cL(S_\tau),\,\cL(S_{\tau_N}))$,
 $A_2=d_{TV}(\cL(S_{\tau_N}),\,\cL(S'_{\tau'_N}))$,
 $A_3=d_{TV}(\cL(S'_{\tau'_N}),\,\cL(S^*_{\tau^*_N}))$ ,
 $A_4=d_{TV}(\cL(S^*_{\tau^*_N}),\,\cL(S^*_{\tau^*}))$ and
 $A_5=d_{TV}(\mbox{Geo}(\varrho),\,\mbox{Geo}(\rho))$.

 Introduce random vectors
 $\bfX_{N,\al}=\{ X_{n,\al},\, 0\leq n\leq N\},\,\al=0,1$, $\bfX_N=
 \{\bfX_{N,0},\,\bfX_{N,1}\}$, $\bfX'_{N,\al}=\{ X'_{n,\al},\, 0\leq n\leq N\},
 \,\al=0,1$, $\bfX'_N=\{\bfX'_{N,0},\,\bfX'_{N,1}\}$, $\bfY_{N,\al}=
 \{ Y_{n,\al},\, 0\leq n\leq N\},\,\al=0,1$ and $\bfY_N=\{\bfY_{N,0},\,
 \bfY_{N,1}\}$. Observe that the event $\{ S_\tau\ne S_{\tau_N}\}$ can
 occur only if $\tau>N$. Also, we can write $\{\tau>N\}=\{ X_{n,0}=0\,\,
 \mbox{for all}\,\, n=0,1,...,N\}$ and $\{\tau'>N\}=\{ X'_{n,0}=0\,\,
 \mbox{for all}\,\, n=0,1,...,N\}$ Hence,
 \begin{eqnarray}\label{5.4}
 &A_1\leq P\{\tau>N\}= P\{\tau'>N\}+|P\{ X_{n,1}=0\,\,\mbox{for}\,\,
 n=0,1,...,N\}\\
 &-P\{ X'_{n,1}=0\,\,\mbox{for}\,\, n=0,1,...,N\}|\leq P\{\tau'>N\}+
 d_{TV}(\cL(\bfX_{N,1}),\,\cL(\bfX'_{N,1}))\nonumber
 \end{eqnarray}
 and similarly,
 \begin{equation}\label{5.5}
 P\{\tau'>N\}\leq P\{\tau^*>N\}+d_{TV}(\cL(\bfX'_{N,1}),\,\cL(\bfY_{N,1})).
 \end{equation}
 Since $Y_{n,1},\, n=0,1,...$ are i.i.d. random variables we obtain that
 \begin{equation}\label{5.6}
 P\{\tau^*>N\}=(P\{ Y_{0,1}=0\})^{N}=(1-(\Phi(\Gam_1))^\ell)^{N}.
 \end{equation}

 Next, we claim that
 \begin{eqnarray}\label{5.7}
 &d_{TV}(\cL(\bfX'_{N,1}),\,\cL(\bfY_{N,1}))\leq d_{TV}(\cL(\bfX'_N),\,
 \cL(\bfY_N))\\
 &\leq\sum_{0\leq n\leq N,\al=0,1}d_{TV}(\cL(X'_{n,\al}),\,\cL(Y_{n,\al})).
 \nonumber \end{eqnarray}
 The first inequality above is clear and the second one holds true since for
 any Borel probability measures $\mu_1,\mu_2$ and $\tilde\mu_1,\tilde\mu_2$
 on Borel measurable spaces $\cX$ and $\tilde\cX$, respectively, (see, for instance, \cite{KR2}),
 \[
 d_{TV}(\mu_1\times\tilde\mu_1,\,\mu_2\times\tilde\mu_2)\leq d_{TV}(\mu_1,\mu_2)+d_{TV}(\tilde\mu_1,\tilde\mu_2).
 \]

 Now, for any $n\in U_{N,M}$ and $\al=0,1$,
 \begin{eqnarray}\label{5.8}
 &d_{TV}(\cL(X'_{n,\al},\,\cL(Y_{n,\al}))=|P\{ X'_{n,\al}=1\}-
 P\{ Y_{n,\al}=1\}|\\
 &=|P\{\xi_{q_{i,N}(n)}\in\Gam_\al\,\,\mbox{for}\,\, i=1,...,\ell\}-
 (\Phi(\Gam_\al))^\ell|\leq\big((1+\psi(M))^\ell-1\big)(\Phi(\Gam_\al))^\ell
 \nonumber\end{eqnarray}
 where the last inequality follows from Lemma 3.2 in \cite{KR1}. By (\ref{3.2}), (\ref{4.5})
 and (\ref{5.8}) for any positive integer $N$ we can write
 \begin{eqnarray}\label{5.9}
 &d_{TV}(\cL(\bfX'_N),\,\cL(\bfY_N))\leq((\Phi(\Gam_0))^\ell+(\Phi(\Gam_1))^\ell))
 N\big((1+\psi(M))^\ell-1\big)\\
 &+KM\ell(\ell-1)(\Phi(\Gam_0)+\Phi(\Gam_1)).\nonumber
 \end{eqnarray}
 Observe that
 \begin{equation}\label{5.10}
 d_{TV}(\cL(\bfX_{N,0}),\cL(\bfX'_{N,0}))\leq d_{TV}(\cL(\bfX_{N}),
 \cL(\bfX'_{N}))\,\,\mbox{and}\,\, A_2\leq d_{TV}(\cL(\bfX_{N}),\cL(\bfX'_{N})).
 \end{equation}
 The first inequality in (\ref{5.10}) is clear and the second one follows
 from the fact that $S_{\tau_N}=f(\bfX_N)$ and $S'_{\tau'_N}=f(\bfX'_N)$
 for a certain function $f:\,\{ 0,1\}^{2L}\to\{ 1,...,N\}$. We will
 estimate next $d_{TV}(\cL(\bfX_N),\,\cL(\bfX'_N))$ relying on \cite{AGG}.

  By Theorem 3 in \cite{AGG},
  \begin{equation}\label{5.11}
  d_{TV}(\cL(\bfX_N),\,\cL(\bfX'_N))\leq 2b_1+2b_2+b_3+2\sum_{0\leq n\leq N,
  \al=0,1}p^2_{n,\al}
  \end{equation}
  where for $\al=0,1$ and $n\in U_{N,M}$,
  \begin{equation}\label{5.12}
  p_{n,\al}=P\{ X_{n,\al}=1\}=P\{\xi_{q_{i,N}(n)}\in\Gam_\al\,\,\mbox{for}\,\,
  i=1,...,\ell\}\leq (1+\psi(M))^\ell (\Phi(\Gam_\al))^\ell
  \end{equation}
  with the latter inequality satisfied by Lemma 3.2 in \cite{KR1}. In order to
  define $b_1,b_2$ and $b_3$ we introduce the set
  \[
  B^{N,R}_{n}=B^{N,R}_{n,N}=\{(l,0),\,(l,1):\, 0\leq l\leq N,\,\del_N(l,n)\leq R\}
  \]
  where an integer $R>0$ is another parameter. Set also
  $I_N=\{(n,\al):\, 0\leq n\leq N,\,\al=0,1\}$. Then
  \begin{equation}\label{5.13}
  b_1=\sum_{(n,\al)\in I_N}\sum_{(l,\be)\in B_{n}^{N,R}}p_{n,\al}
  p_{l,\be},
  \end{equation}
  \begin{equation}\label{5.14}
  b_2=\sum_{(n,\al)\in I_N}\sum_{(n,\al)\ne(l,\be)\in B_{n}^{N,R}}
  p_{(n,\al),(l,\be)},
  \end{equation}
 where $p_{(n,\al),(l,\be)}=E(X_{n,\al}X_{l,\be})$, and
 \begin{equation}\label{5.15}
 b_3=\sum_{(n,\al)\in I_N}s_{n,\al}
 \end{equation}
 where
 \[
 s_{n,\al}=E\big\vert E\big(X_{n,\al}-p_{n,\al}|\sig\{ X_{l,\be}:\,(l,\be)\in
 I_N\setminus B^{N,R}_{n}\}\big)\big\vert.
 \]
 By Assumption \ref{ass2.1}, for any $i,j,n$ and $k$ there exists at most $K$ of $l$'s
 such that $q_{i,N}(n)-q_{j,N}(l)=k$. It follows from here that
  \begin{equation}\label{5.16}
  \# B^{N,R}_{n}\leq K\ell^2R.
  \end{equation}
   It follows from (\ref{3.2}), (\ref{5.12}), (\ref{5.13}) and (\ref{5.16}), similarly to (\ref{4.6}), that
  \begin{eqnarray}\label{5.17}
  &b_1\leq 2N\ell^2R(1+\psi(1))^{2\ell}((\Phi(\Gam_0))^{2\ell}+(\Phi(\Gam_1))^{2\ell})\\
  &+K\ell^2(1+\psi(1))^{\ell+1}(N+R)((\Phi(\Gam_0))^{\ell+1}+(\Phi(\Gam_1))^{\ell+1})\nonumber\\
  &+K^2\ell^4(1+\psi(1))^2((\Phi(\Gam_0))^{2}+(\Phi(\Gam_1))^{2}).
  \nonumber\end{eqnarray}

  Since $\Gam_0\cap\Gam_1=\emptyset$,
  \begin{eqnarray}\label{5.18}
  &p_{(n,\al),(l,\be)}=P\{ X_{n,\al}=X_{l,\be}=1\}=0\,\,\mbox{if}\,\, n=l,\,\be=1-\al\\
  &\mbox{and always}\,\, p_{(n,\al),(l,\be)}\leq p_{(l,\be)}.\nonumber\\
  \end{eqnarray}
  Similarly to Section \ref{sec4} we estimate $p_{(n,\al),(l,\be)}$ by the right hand sides of (\ref{4.7})--(\ref{4.11})
  in the corresponding cases replacing  $\Phi(\Gam)$ there by $\Phi(\Gam_0)+\Phi(\Gam_1)$ here, namely,
   if $\del_N(n,l)=k\geq 1$ and $n,l\in U_{N,M}$ we estimate it via the right hand side of (\ref{4.7}),
   if either $n\in U_{N,M}$ or $l\in U_{N,M}$ and either $\del_N(n,l)=k\geq 1$ or $\del_N(n,l)=0$ and $\sig_N(n,l)\geq 1$
   we estimate it via the right hand side of (\ref{4.8}),
  if $\del_N(n,l)=k\geq 1$ and $n,l\not\in U_{N,M}$ we estimate it via the right hand side of (\ref{4.9}),
  if $\sig_N(n,l)=0$ and either $n\in U_{N,M}$ or $l\in U_{N,M}$ we estimate it by the right hand side of (\ref{4.10})
  and, finally, if $\sig_N(n,l)=0$ and $n,l\not\in U_{N,M}$ we estimate it via the right hand side of (\ref{4.11}). These estimates
  together with counting estimates of Section \ref{sec3} yield
  \begin{eqnarray}\label{5.19}
  &b_2\leq 2NKR\ell^2(1+\psi(1))^{2\ell}(\Phi(\Gam_0)+\Phi(\Gam_1))^{2\ell}\\
&+2KM\ell^2(N+KR\ell^2)(1+\psi(1))^{\ell+1}(\Phi(\Gam_0)+\Phi(\Gam_1))^{\ell+1}\nonumber\\
&+K^2M^2\ell^4(\Phi(\Gam_0)+\Phi(\Gam_1))+K^2M^2\ell^4(1+\psi(1))^2(\Phi(\Gam_0)+\Phi(\Gam_1))^2\nonumber\\
&+2K^2MR\ell^4(1+\psi(M))^\ell(\Phi(\Gam_0)+\Phi(\Gam_1))^\ell.\nonumber
 \end{eqnarray}

 In the same way as in Section \ref{sec4} we obtain that
  \begin{eqnarray}\label{5.20}
  &s_{n,\al}\leq 2^{2(\ell+2)}(2-(1+\psi(R))^{\ell+1})^{-2}\psi(R)E|X_{n,\al}
  -p_{n,\al}|\\
  &\leq 2^{2\ell+5}(2-(1+\psi(R))^{\ell+1})^{-2}\psi(R)p_{n,\al}
  \nonumber\end{eqnarray}
  where $s_{n,\al}$ is the same as in (\ref{5.15}).
  Hence, by (\ref{5.12}), (\ref{5.15}) and (\ref{5.20}),
  \begin{eqnarray}\label{5.21}
  &b_3\leq 2^{2\ell+5}(2-(1+\psi(R))^{\ell+1})^{-2}\psi(R)\big( 2KM\ell(\ell-1)(\Phi(\Gam_0)+\Phi(\Gam_1))\\
  &+N(1+\psi(M))^\ell((\Phi(\Gam_0))^\ell+(\Phi(\Gam_1))^\ell)\big).\nonumber
  \end{eqnarray}

  Next, in the same way as in the estimate of $A_2$ we conclude that
  \begin{equation}\label{5.22}
  A_3\leq d_{TV}(\cL(\bfX'_N),\cL(\bfY_N))
  \end{equation}
  which together with (\ref{5.9}) estimates $A_3$.

  As in the estimate of $A_1$ we see that
  \begin{equation}\label{5.23}
  A_4\leq P\{\tau^*>N\}\leq (1-(\Phi(\Gam_1))^\ell)^{N}
  \end{equation}
  since $Y_{n,0},\, n=0,1,...$ are i.i.d. random variables.

  Since $\varrho>\rho$ we obtain
  \begin{eqnarray}\label{5.24}
  &A_5\leq\sum_{k=0}^\infty |\varrho(1-\varrho)^k-\rho(1-\rho)^k|\leq
  2\sum_{k=1}^\infty((1-\rho)^k-(1-\varrho)^k)\\
  &=2(1-\rho)\rho^{-1}-2(1-\varrho)\varrho^{-1}=
  2\frac {\varrho-\rho}{\rho\varrho}=2(\Phi(\Gam_0))^\ell. \nonumber
  \end{eqnarray}
 Collecting (\ref{5.3})--(\ref{5.12}), (\ref{5.17}), (\ref{5.20}),
 (\ref{5.19}) and (\ref{5.20})--(\ref{5.24}) we derive (\ref{2.8}).
  \qed

 In order to prove Corollary \ref{cor2.5} we rely on the estimate (\ref{2.8})
 with $\Gam=\Gam_N$ and $\Del=\Del_N$ choosing  $M=M_N\to\infty$ and
 $R=R_N\to\infty$ as $N\to\infty$ so that
 \begin{eqnarray}\label{5.25}
 &\lim_{N\to\infty}N(\Phi(\Del_N))^\ell=\infty,\\
 &\lim_{N\to\infty}\big((\Phi(\Gam_N))^\ell+(\Phi(\Del_N))^\ell\big)N\psi(M_N)=0,\nonumber\\
 &\lim_{N\to\infty}N\psi(R_N)\big((\Phi(\Gam_N))^\ell+(\Phi(\Del_N))^\ell\big)=0,\nonumber\\
 &\lim_{N\to\infty}N(R_N+M_N)\big((\Phi(\Gam_N))^{\ell+1}+(\Phi(\Del_N))^{\ell+1}\big)=0\quad\mbox{and}\nonumber\\
 &\lim_{N\to\infty}\big(M^2_N(\Phi(\Gam_N)+\Phi(\Del_N))+M_NR_N\big((\Phi(\Gam_N))^\ell+(\Phi(\Del_N))^\ell\big)\big)=0.\nonumber
 \end{eqnarray}
 which is clearly possible since $\psi(n)\to 0$ as $n\to\infty$. This together with (\ref{2.8})
  yields (\ref{2.10}).
 \qed

\section{Poisson distribution limits for shifts}\label{sec6}\setcounter{equation}{0}

Let $V\in\cF_{0,n-1}$ and set $X_{k,N}=X^V_{k,N}=\prod_{i=1}^\ell\bbI_V\circ T^{q_{i,N}(k)}$.
Then $S_N=S_N^V=\sum_{k=1}^NX_{k,N}$. Set $p_{k,N}=P\{ X_{k,N}=1\}$ and $p_{k,l,N}=P\{ X_{k,N}=1$ and $X_{l,N}=1\}$.
 Then, again, by Theorem 1 from \cite{AGG} we obtain
\begin{equation}\label{6.1}
d_{TV}(\cL(S_N),\, Pois(\la_N))\leq b_1+b_2+b_3
\end{equation}
where $b_1,\, b_2$ and $b_3$ are defined by (\ref{4.2}) and (\ref{4.3}) though their estimates will proceed
now somewhat differently than in Section \ref{sec4}.

Since $V\in\cF_{0,n-1}$, it follows from Lemma 3.2 in \cite{KR1} that for any $k\in U_{N,n}$,
\begin{equation}\label{6.2}
p_{k,N}=EX_{k,N}\leq (1+\psi(1))^\ell(P(V))^\ell
\end{equation}
while when $k\not\in U_{N,n}$ we can always write
\begin{equation}\label{6.3}
p_{k,N}\leq E(\bbI_V\circ T^{q_{1,N}(k)})=P(V).
\end{equation}
Hence, by (\ref{3.2}) and (\ref{3.4}) we conclude that
\begin{equation}\label{6.4}
b_1=\sum_{k=1}^N\sum_{l\in B^R_{k,N}}p_{k,N}p_{l,N}\leq K^2R\ell^4P(V)+NKR\ell^2(P(V))^{\ell+1}.
\end{equation}

In order to estimate $p_{k,l,N}$ we observe that if $|i-j|<\pi(V)$ then $(\bbI_V\circ T^i)(\bbI_V\circ T^j)=0$. Hence,
$p_{k,l,N}=0$ if $\del_N(k,l)<\pi(V)$.
Now suppose that $\del_N(k,l)=d$ with $\pi(V)\leq d<n$,
\[
q_{i_1,N}(k)\leq q_{i_2,N}(k)\leq...\leq q_{i_\ell,N}(k)\,\,\mbox{and}\,\, q_{j_1,N}(l)\leq q_{j_2,N}(l)\leq...\leq q_{j_\ell,N}(l).
\]
Assume that the pair $k,l$ does not belong to the exceptional set$D_N$ of cardinality at most $K$ appearing in Assumption \ref{ass2.1}(iii).
Since $\del_N(k,l)=d\geq\pi(V)$, it follows that
\begin{equation}\label{6.5}
\mbox{either}\,\,\, q_{j_1,N}(l)\leq q_{i_1,N}(k)-d\,\,\,\mbox{or}\,\,\, q_{j_\ell,N}(l)\geq q_{i_\ell,N}(k)+d
\end{equation}
and in view of Assumption \ref{ass2.1}(iii) only one of these inequalities can hold true. Assume, for instance, that the first inequality
in (\ref{6.5}) holds true and let $r=q_{i_1,N}(k)-q_{j_1,N}(l)$.
Then $r\geq d\geq\pi(V)$. If $r\geq n$ then by Lemma 3.2 from \cite{KR1} (essentially, by the definition of the $\psi$-mixing coefficient itself),
\begin{equation}\label{6.6}
p_{k,l,N}=E(X_{k,N}X_{l,N})\leq E(X_{k,n}\bbI_V\circ T^{q_{j_1,N}(l)})\leq (1+\psi(1))p_{k,N}P(V).
\end{equation}

Suppose that $\pi(V)\leq r<n$. Then the sets $Q_0=\{ q_{j_1,N}(l),\, q_{j_1,N}(l)+1,...,q_{j_1,N}(l)+n-1\}$ and
$Q_1=\{ q_{i_1,N}(k)+n-r,\, q_{i_1,N}(k)+n-r+1,...,q_{i_1,N}(k)+n-1\}$ are disjoint, and so it follows by Lemma 3.2
from \cite{KR1} that in this case,
\begin{eqnarray}\label{6.7}
&p_{k,l,N}=E(X_{k,N}X_{l,N})\leq E(\bbI_V\circ T^{q_{j_1,N}(l)}\bbI_V\circ T^{q_{i_1,N}(k)})\\
&\leq E(\bbI_V\circ T^{q_{j_1,N}(l)}\bbI_{T^{n-r}V}\circ T^{q_{i_1,N}(k)+n-r}V)\leq (1+\psi(1))P(V)P(T^{n-r}V)\nonumber
\end{eqnarray}
where we used that $V\subset T^{-a}T^aV$ for any integer $a\geq 0$. If, in addition, $k\in U_{N,n}$ then the sets
$Q_0=\{ q_{j_1,N}(l),\, q_{j_1,N}(l)+1,...,q_{j_1,N}(l)+n-1\}$,
$Q_1=\{ q_{i_1,N}(k)+n-r,\, q_{i_1,N}(k)+n-r+1,...,q_{i_1,N}(k)+n-1\}$ and
$Q_m=\{ q_{i_m,N}(k),q_{i_m,N}(k)+1,...,q_{i_m,N}(k)+n-1\}$ , $m=1,...,\ell$ are disjoint and we obtain then from Lemma 3.2
 in \cite{KR1} that
\begin{equation}\label{6.8}
p_{k,l,N}=E(X_{k,N}X_{l,N})\leq (1+\psi(1))^\ell (P(V))^\ell P(T^{n-r}V).
\end{equation}

 If the second inequality in (\ref{6.5}) holds true then we obtain (\ref{6.6}) if $r=q_{j_\ell,N}(l)-q_{i_\ell,N}(k)\geq n$,
  while if $\pi(V)\leq r<n$ then we use that the sets $Q_\ell=\{ q_{i_\ell,N}(k),q_{i_\ell,N}(k)+1,...,q_{i_\ell,N}(k)+n-1\}$ and
  $Q_{\ell+1}=\{ q_{j_\ell,N}(l)+n-r,q_{j_\ell,N}(l)+n-r+1,...,q_{j_\ell,N}(l)+n-1\}$ are disjoint which yields (\ref{6.7}) by
  Lemma 3.2 from \cite{KR1}. If, in addition, $k\in U_{N,n}$ then all sets  $Q_m=\{ q_{i_m,N}(k),q_{i_m,N}(k)+1,...,q_{i_m,N}(k)+n-1\}$,
  $m=1,...,\ell$ and $Q_{\ell+1}=\{ q_{j_\ell,N}(l)+n-r,q_{j_\ell,N}(l)+n-r+1,...,q_{j_\ell,N}(l)+n-1$ are disjoint, and so by Lemma 3.2
  from \cite{KR1} we obtain the estimate (\ref{6.8}) again. Finally, suppose that $\del_N(k,l)=d\geq n$. Then, applying Lemma 3.2 from
  \cite{KR1} we see that the estimate (\ref{6.6}) holds true again. Observe that by Assumption \ref{ass2.1}(i) for any $N\geq 1$, $i=1,...,\ell$
  and integers $k\geq 0$ and $r$,
  \begin{equation}\label{6.9}
  \#\{ l\geq 0:\, q_{i,N}(k)-q_{j,N}(l)=r\,\,\mbox{for some}\,\, 1\leq j\leq\ell\}\leq\ell K.
  \end{equation}
  Now, it follows from (\ref{3.2})--(\ref{3.4}), (\ref{6.2}), (\ref{6.3})
  and (\ref{6.6})--(\ref{6.9}) that
  \begin{eqnarray}\label{6.10}
&\quad b_2=\sum_{k=1}^N\sum_{k\ne l\in B_{k,N}^R}p_{k,l,N}=\sum_{1\leq k\leq N,k\not\in U_{N,n}}\sum_{l:\,\pi(V)\leq\del_N(k,l)<R}p_{k,l,N}\\
&+\sum_{k\in U_{N,n}}\sum_{l:\,\pi(V)\leq\del_N(k,l)<R}p_{k,l,N}\leq K^2\ell^4(1+\psi(1))nR(P(V))^2\nonumber\\
&+K^2\ell^3(1+\psi(1))nP(V)\sum_{r=\pi(V)}^{n-1}P(T^{n-r}V)+K\ell^2(1+\psi(1))^\ell NR(P(V))^{\ell+1}\nonumber\\
&+K\ell((1+\psi(1))^\ell N(P(V))^\ell\big)\sum_{r=\pi(V)}^{n-1}P(T^{n-r}V)=\big( K^2\ell^3(1+\psi(1))nP(V)\nonumber\\
&+K\ell((1+\psi(1))^\ell N(P(V))^\ell\big)\big(R\ell P(V)+\sum_{r=\pi(V)}^{n-1}P(T^{n-r}V)\big).
\nonumber\end{eqnarray}

Next, we estimate $s_{k,N}$ and $b_3$ defined by (\ref{4.3}). Let $\cG=\cG_k$ be the $\sig$-algebra generated by the
sets $T^{-q_{i,N}(l)}V,\, i=1,...,\ell;\, l\in\{ 1,...,N\}\setminus B^R_{k,N}$ and $\cH=\cH_k$ be the $\sig$-algebra generated
by the sets $T^{-q_{i,N}(k)}V,\, i=1,...,\ell$. Since $\del_N(k,l)\geq R$ for $l\not\in B^R_{k,N}$ and $V\in\cF_{0,n-1}$,
we derive from Lemma 3.3 in \cite{KR1} that for $n<R<N$,
\begin{equation}\label{6.11}
\psi(\cG,\cH)\leq 2^{2(\ell+2)}\psi(R-n)(2-(1+\psi(R-n))^{\ell+1})^{-2},
\end{equation}
and so
\begin{eqnarray}\label{6.12}
&s_{k,N}\leq 2^{2(\ell+2)}\psi(R-n)(2-(1+\psi(R-n))^{\ell+1})^{-2}E|X_{k,N}-p_{k,N}|\\
&\leq 2^{2\ell+5}\psi(R-n)p_{k,N}(2-(1+\psi(R-n))^{\ell+1})^{-2}.\nonumber
\end{eqnarray}
Hence, by (\ref{3.2}), (\ref{6.2}), (\ref{6.3}) and (\ref{6.12}),
\begin{eqnarray}\label{6.13}
&b_3=\sum_{k=1}^Ns_{k,N}\leq 2^{2\ell+5}\psi(R-n)(2-(1+\psi(R-n))^{\ell+1})^{-2}\\
&\times\big( K\ell^2nP(V)+(1+\psi(1))^\ell N(P(V))^\ell\big)\nonumber
\end{eqnarray}
Finally, collecting (\ref{6.1}), (\ref{6.4}), (\ref{6.9}) and (\ref{6.13}) we derive (\ref{2.13}) completing the proof of
Theorem \ref{thm2.6}.
\qed

Corollary \ref{cor2.7} follows from the estimate (\ref{2.13}) choosing $R=R_L=2n_L$ and in view of (\ref{2.12}) we obtain Corollary \ref{cor2.7}
for $V_{N_L}=A^\eta_{n_L}$ since
\begin{equation}\label{6.14}
\pi(A^\eta_n)\to\infty\quad\mbox{as}\quad n\to\infty
\end{equation}
when $\eta$ is a nonperiodic sequence. Indeed, $\pi(A^\eta_n)$ is, clearly, nondecreasing in $n$, and so $\lim_{n\to\infty}\pi(A^\eta_n)=r$ exists.
If $r<\infty$ then there exists $n_0\geq 1$ such that $\pi(A^\eta_n)=r$ for all $n\geq n_0$ which means that $\eta$ is periodic with the period $r$,
and so $r=\infty$ since $\eta$ is not periodic.
\qed

\section{Geometric distribution limits for shifts}\label{sec7}\setcounter{equation}{0}

It will be convenient to set $V^{(0)}=V\in\cF_{0,n-1}$, $V^{(1)}=W\in\cF_{0,m-1}$ and
\[
X^{(\al)}_{k,N}=\prod_{i=1}^\ell\bbI_{V^{(\al)}}\circ T^{q_{i,N}(k)},\,\,\al=0,1
\]
 so that
 \[
 \tau=\tau_{V^{(1)}}=\min\{ k\geq 1:\, X^{(1)}_{k,N}=1\}\,\,\,\mbox{and}\,\,\,\Sig_N^{V^{(0)},V^{(1)}}=\sum_{k=1}^\tau X^{(0)}_{k,N}.
 \]
 Set also $S_L=\sum_{k=1}^LX_{k,N}^{(0)}$, so that $S_\tau=\Sig_N^{V^{(0)},V^{(1)}}$, and denote $\tau_N=\min(\tau,N)$. Let $\{ Y^{(\al)}_{k,N}:\,
 k\geq 1,\,\al=0,1\}$ be a sequence of independent Bernoulli random variables such that $Y^{(\al)}_{k,N}$ has the same distribution
 as $X^{(\al)}_{k,N}$. Set
 \[
 S'_L=\sum_{k=1}^LY^{(0)}_{k,N},\,\,\tau'=\min\{ k\geq 1:\, Y^{(1)}_{k,N}=1\}\,\,\mbox{and}\,\,\tau'_N=\min(\tau',N).
 \]
 Let now $\{ Z^{(\al)}_{k,N}:\, k\geq 1\}$, $\al=0,1$ be two independent of each other sequences of i.i.d. Bernoulli random variables
 such that
 \begin{equation}\label{7.1}
 P\{ Z^{(\al)}_{k,N}=1\}=(P(V^{\al)}))^\ell=1-P\{ Z^{(\al)}_{k,N}=0\},\,\,\al=0,1.
 \end{equation}
 We can and will assume that all above random variables are defined on the same (sufficiently large) probability space. Set also
 \[
 S^*_L=\sum_{k=1}^L Z^{(0)}_{k,N},\,\,\tau'=\min\{ k\geq 1:\, Z^{(1)}_{k,N}=1\}\,\,\mbox{and}\,\,\tau^*_N=\min(\tau^*,N).
 \]
 By Lemma 3.1 from \cite{KR1} the sum $S^*_{\tau^*}$ has the geometric distribution with the parameter
 \begin{equation}\label{7.2}
 \vrho=\frac {(P(V^{(1)}))^\ell}{(P(V^{(1)}))^\ell+(P(V^{(0)}))^\ell(1-(P(V^{(1)}))^\ell)}>\rho
 \end{equation}
 where $\rho=(P(V^{(1)}))^\ell\big(P(V^{(1)}))^\ell+(P(V^{(0)}))^\ell\big)^{-1}$.

 Next, we can write
\begin{equation}\label{7.3}
d_{TV}(\cL(S_\tau),\,\mbox{Geo}(\rho))\leq A_1+A_2+A_3+A_4+A_5
\end{equation}
where $A_1=d_{TV}(\cL(S_\tau),\,\cL(S_{\tau_N}))$,
 $A_2=d_{TV}(\cL(S_{\tau_N}),\,\cL(S'_{\tau'_N}))$,
 $A_3=d_{TV}(\cL(S'_{\tau'_N}),\,\cL(S^*_{\tau^*_N}))$ ,
 $A_4=d_{TV}(\cL(S^*_{\tau^*_N}),\,\cL(S^*_{\tau^*}))$ and
 $A_5=d_{TV}(\mbox{Geo}(\varrho),\,\mbox{Geo}(\rho))$.

 Introduce random vectors
 $\bfX_{N}^{(\al)}=\{ X_{k,n}^{(\al)},\, 1\leq k\leq N\},\,\al=0,1$, $\bfX_N=
 \{\bfX_{N}^{(0)},\,\bfX_{N}^{(1)}\}$, $\bfY_{N}^{(\al)}=\{ Y_{n,N}^{(\al)},\, 1\leq k\leq N\},
 \,\al=0,1$, $\bfY_N=\{\bfY_{N}^{(0)},\,\bfY_{N}^{(1)}\}$, $\bfY_{N}^{(\al)}=
 \{ Y_{n,N}^{(\al)},\, 1\leq k\leq N\},\,\al=0,1$ and $\bfY_N=\{\bfY_{N}^{(0)},\,
 \bfY_{N}^{(1)}\}$. Observe that the event $\{ S_\tau\ne S_{\tau_N}\}$ can
 occur only if $\tau>N$. Also, we can write $\{\tau>N\}=\{ X_{n,N}^{(1)}=0\,\,
 \mbox{for all}\,\, k=1,...,N\}$ and $\{\tau'>N\}=\{ Y_{n,0}^{(1)}=0\,\,
 \mbox{for all}\,\, k=1,...,N\}$ Hence,
 \begin{eqnarray}\label{7.4}
 &A_1\leq P\{\tau>N\}= P\{\tau'>N\}+|P\{ X_{n,N}^{(1)}=0\,\,\mbox{for}\,\,
 n=1,...,N\}\\
 &-P\{ Y_{n,N}^{(1)}=0\,\,\mbox{for}\,\, n=0,1,...,N\}|\leq P\{\tau^*>N\}\nonumber\\
 &+ d_{TV}(\cL(\bfY_{N}),\,\cL(\bfZ_{N}))+d_{TV}(\cL(\bfX_{N}),\,\cL(\bfY_{N})).\nonumber
 \end{eqnarray}

 Since $Z^{(1)}_{k,N},\, k=0,1,...$ are i.i.d. random variables we obtain that
 \begin{equation}\label{7.5}
 P\{\tau^*>N\}=(1-(P(V^{(1)})^\ell))^{N}.
 \end{equation}
It is also not difficult to understand (see p.p. 1534--1535 in \cite{KR1}) that
\begin{equation}\label{7.6}
d_{TV}(\cL(\bfY_{N}),\,\cL(\bfZ_{N}))\leq\sum_{1\leq k\leq N,\,\al=0,1}d_{TV}(\cL(Y_{k,N}^{(\al)}),\cL(Z^{(\al)}_{k,N})).
\end{equation}
If $k\in U_{N,n\vee m}$ then by (\ref{7.6}) and Lemma 3.2 from \cite{KR1} similarly to (\ref{5.9}) we obtain that,
\begin{eqnarray}\label{7.7}
&d_{TV}(\cL(Y^{(\al)}_{k,N}),\,\cL(Z^{(\al)}_{k,N}))=|P\{ Y^{(\al)}_{k,N}=1\}-P\{ Z^{(\al)}_{k,N}=1\}|\\
&=|P\{\cap_{i=1}^\ell T^{-q_{i,N}(k)}V^{(\al)}\}-(P(V^{(\al)}))^\ell|\leq((1+\psi(n\vee m))^\ell-1)(P(V^{(\al)}))^\ell.
\nonumber\end{eqnarray}
It follows from (\ref{3.2}), (\ref{7.6}) and (\ref{7.7}) that
\begin{eqnarray}\label{7.8}
&d_{TV}(\cL(\bfY_{N},\,\cL(\bfZ_{N}))\leq K\ell^2(n\vee m)(P(V^{(0}))+P(V^{(1)})\\
&+N\big((P(V^{(0})))^\ell+(P(V^{(1)}))^\ell\big)((1+\psi(n\vee m))^\ell-1).
\nonumber\end{eqnarray}

Next, we observe that by Theorem 3 in \cite{AGG},
\begin{equation}\label{7.9}
A_2\leq d_{TV}(\cL(\bfX_{N}),\,\cL(\bfY_{N}))\leq 2b_1+2b_2+b_3+2\sum_{1\leq k\leq N,\,\al=0,1}(p_{k,N}^{(\al)})^2
\end{equation}
where $p_{k,N}^{(\al)}=P\{ X_{k,N}^{(\al)}=1\}$ and if $k\in U_{N,n\vee m}$ then (in the same way as in (\ref{6.2}) by Lemma 3.2 in \cite{KR1},
\begin{equation}\label{7.10}
p_{k,N}^{(\al)}\leq (1+\psi(n\vee m))^\ell(P(V^{(\al)}))^\ell,
\end{equation}
while the definitions of $b_1,b_2$ and $b_3$ are similar to Section \ref{sec6} taking into account the additional parameter $\al$.
Namely, setting
\[
B^{N,R}_k=\{ (l,0),\, (l,1):\, 1\leq l\leq N,\,\del(k,l)\leq R\},\, p^{\al,\be}_{k,l,N}=E(X^{(\al)}_{k,N}X^{(\be)}_{l,N})
\]
and $I_N=\{(k,\al):\, 1\leq k\leq N,\,\al=0,1\}$ we have
\begin{equation}\label{7.11}
b_1=\sum_{(k,\al)\in I_N}\sum_{(l,\be)\in B^{N,R}_k}p^{(\al)}_{k,N}p^{(\be)}_{l,N},
\end{equation}
\begin{equation}\label{7.12}
b_2=\sum_{(k,\al)\in I_N}\sum_{(k,\al)\ne(l,\be)\in B^{N,R}_k}p^{(\al,\be)}_{k,l,N}\quad\mbox{and}
\end{equation}
\begin{equation}\label{7.13}
b_3=\sum_{(k,\al)\in I_N}s^{(\al)}_{k,N}\,\,\,\mbox{where}
\end{equation}
\[
s^{(\al)}_{k,N}=E\big\vert E\big(X^{(\al)}_{k,N}-p^{(\al)}_{k,N}|
\sig\{ X^{(\be)}_{l,N}:\, (l,\be)\in I_N\setminus B^{N,R}_k\}\big)\big\vert.
\]

Since we always have
\begin{equation}\label{7.14}
p^{(\al)}_{k,N}\leq P(V^{(\al)})
\end{equation}
and (\ref{7.10}) holds true when $k\in U_{N,n\vee n}$, it follows taking into account (\ref{3.2}) and (\ref{3.4}) that
\begin{eqnarray}\label{7.15}
&b_1\leq K\ell^2(1+\psi(n\vee m))^\ell RN\big((P(V^{(0}))^\ell+(P(V^{(1)}))^\ell\big)\\
&\times(P(V^{(0})+P(V^{(1)}))+K^2\ell^2(n\vee m)^2(P(V^{(0})+P(V^{(1)}))^2.\nonumber
\end{eqnarray}
In order to estimate $p^{\al,\be}_{k,l,N}$ (and, eventually, $b_2$) we will essentially repeat the arguments from Section \ref{sec6}.
First, observe that
\[
(\bbI_{V^{(0)}}\circ T^i)(\bbI_{V^{(1)}}\circ T^j)=0\,\,\mbox{if}\,\, |i-j|<\ka_{V^{(0)},V^{(1)}}.
\]
Hence, $p^{\al,\be}_{k,l,N}=0$ if $\del_N(k,l)<\ka_{V^{(0)},V^{(1)}}$. Now suppose that
$\del_N(k,l)=d\geq\ka_{V^{(0)},V^{(1)}}$, $q_{i_1,N}(k)\leq q_{i_2,N}(k)\leq...\leq q_{i_\ell,N}(k)$
and $q_{j_1,N}(l)\leq q_{j_2,N}(l)\leq...\leq q_{j_\ell,N}(l)$. Assume that the pair $k,l$ does not belong to the exceptional set $D_N$
of cardinality at most $K$ appearing in Assumption \ref{ass2.1}(iii). Then we have to deal with two alternatives from (\ref{6.5}).

If the first inequality in (\ref{6.5}) holds true and $r=q_{i_1,N}(k)-q_{j_1,N}(l)\geq n\vee m$ then by Lemma 3.2 from \cite{KR1},
\begin{equation}\label{7.16}
p^{\al,\be}_{k,l,N}=E(X^{(\al)}_{k,N}X^{(\be)}_{l,N})\leq (1+\psi(1))p^{(\al)}_{k,N}P(V^{(\be)}).
\end{equation}
If, on the other hand, $\ka_{V^{(0)},V^{(1)}}\leq r<n\vee m$ then in the same way as in Section \ref{sec6} we obtain that
\begin{equation}\label{7.17}
p^{\al,\be}_{k,l,N}=E(X^{(\al)}_{k,N}X^{(\be)}_{l,N})\leq (1+\psi(1))P(V^{(\al)})P(T^{n\vee m-r}V^{(\be)}).
\end{equation}
If, in addition, $k\in U_{N,n\vee m}$ then in the same way as in (\ref{6.8}),
\begin{equation}\label{7.18}
p^{\al,\be}_{k,l,N}\leq (1+\psi(1))^\ell(P(V^{(\al)}))^\ell P(T^{n\vee m-r}V^{(\be)}).
\end{equation}
If the second inequality in (\ref{6.5}) holds true then we obtain (\ref{7.16}) if $r=q_{j_\ell,N}(l)-q_{i_\ell,N}(k)\geq n\vee m$
while (\ref{7.17}) follows if $\ka_{V^{(0)},V^{(1)}}\leq r<n\vee m$ and if, in addition, $k\in U_{N,n\vee m}$ then we obtain (\ref{7.18}).
Relying on (\ref{3.2})--(\ref{3.4}), (\ref{6.9}), (\ref{7.10}), (\ref{7.12}), (\ref{7.14}) and (\ref{7.16})--(\ref{7.18}) we conclude similarly
to (\ref{6.10}) that
\begin{eqnarray}\label{7.19}
&b_2\leq \bigg(K^2\ell^3(1+\psi(1))(n\vee m)(P(V^{(0)})+P(V^{(1)}))\\
&+K\ell(1+\psi(1))^\ell N\big((P(V^{(0)}))^\ell+(P(V^{(1)}))^\ell\big)\bigg)\bigg(\ell R(P(V^{(0)})+P(V^{(1)}))\nonumber\\
&+\sum_{r=\ka_{V^{(0)},V^{(1)}}}^{n\vee m-1}\big(P(T^{n\vee m-r}V^{(0)})+P(T^{n\vee m-r}V^{(1)})\big)\bigg).
\nonumber\end{eqnarray}
Similarly to (\ref{6.13}) we obtain also that
\begin{eqnarray}\label{7.20}
&b_3\leq 2^{2\ell+5}\psi(R-n\vee m)(2-(1+\psi(R-n\vee m))^{\ell+1})^{-2}\\
&\times\bigg( K\ell^2(n\vee m)(P(V^{(0)})+P(V^{(1)}))\nonumber\\
&+(1+\psi(1))^\ell N\big((P(V^{(0)}))^\ell+(P(V^{(0)}))^\ell\big)\bigg).\nonumber
\end{eqnarray}
These provide the estimate of $A_2$ by (\ref{7.9}), (\ref{7.14}), (\ref{7.15}), (\ref{7.19}) and (\ref{7.20}).

In order to estimate $A_3$ observe that
\begin{eqnarray}\label{7.21}
&A_3\leq d_{TV}(\cL(\bfY_{N}),\,\cL(\bfZ_{N}))\leq\sum_{1\leq k\leq N,\,\al=0,1}d_{TV}(\cL(Y_{k,N}^{(\al)}),\cL(Z^{(\al)}_{k,N}))\\
&=\sum_{1\leq k\leq N,\,\al=0,1}|P\{ Y^{(\al)}_{k,N}=1\}-P\{ Z^{(\al)}_{k,N}=1\}|\nonumber\\
&\leq(\sum_{k\in U_{N,n\vee m},\,\al=0,1}+\sum_{1\leq k\leq N,\, k\not\in U_{N,n\vee m},\,\al=0,1})\big\vert
P\{ Y^{(\al)}_{k,N}=1\}-(P(V^{(\al)}))^\ell\big\vert\nonumber\\
&\leq((1+\psi(n\vee m))^\ell-1)N\big((P(V^{(0)}))^\ell+ (P(V^{(1)}))^\ell\big)\nonumber\\
&+K\ell^2(n\vee m)(P(V^{(0)})+P{(1)}))\nonumber
\end{eqnarray}
where in the last inequality we relied on Lemma 3.2 from \cite{KR1} and on (\ref{3.2}) above. The estimate of $A_4$ we obtain from (\ref{7.5}),
\begin{equation}\label{7.22}
A_4\leq P\{\tau^*>N\}=(1-(P(V^{(1)}))^\ell)^N
\end{equation}
since we are dealing here with an i.i.d. sequence of Bernoulli random variables. The estimate of $A_5$ is the same as in (\ref{5.24}),
\begin{equation}\label{7.23}
A_5\leq 2P(V^{(0)}).
\end{equation}
Finally, combining (\ref{7.3})--(\ref{7.5}), (\ref{7.8})--(\ref{7.15}) and (\ref{7.19})--(\ref{7.23}) we derive (\ref{2.15}) completing the
 proof of Theorem \ref{thm2.8}.       \qed

 Corollary \ref{cor2.9} follows from the estimate (\ref{2.15}) choosing $R=R_L=2(n_L\vee m_L)$ and if $V_L=A_{n_L}^\xi$ and $W_L=A^\eta_{m_L}$ it remains
 only to verify the assertion that $\ka_{A_n^\xi,A_m^\eta}\to\infty$ as $n,m\to\infty$ provided that $\xi,\eta\in\Om_P$ are not periodic and not
 shifts of each other. Indeed, $\pi(A^\xi_n)$, $\pi(A^\eta_m)$ and $\pi(A^\xi_n,A^\eta_m)$ are nondecreasing in $n$ and $m$, and so does
 $\pi(A^\xi_n,A^\eta_m)$. Hence, the limit $r=\lim_{n,m\to\infty}\ka_{A_n^\xi,A_m^\eta}$ exists. If $r<\infty$ then, at least, one of the limits
 $r_1=\lim_{n\to\infty}\pi(A^\xi_n)$, $r_2=\lim_{m\to\infty}\pi(A^\eta_m)$ or $r_3=\lim_{n,m\to\infty}\pi(A^\xi_n,A^\eta_m)$ is finite. If $r_1<\infty$
 then $\xi$ is periodic with the period $r_1$, if $r_2<\infty$ then $\eta$ is periodic with the period $r_2$ and if $r_3<\infty$ then either $T^{r_3}\xi
 =\eta$ or $T^{r_3}\eta=\xi$. Finally, it follows from Lemma 4.1 from \cite{KR2} that (\ref{2.24}) holds true for $P\times P$-almost all $(\xi,\eta)$,
   completing the proof.    \qed


\end{document}